\documentclass{article}

\usepackage{amssymb}
\usepackage{amsmath,enumerate,rotate}
\usepackage{amssymb,latexsym}

\usepackage[english]{babel}

\newtheorem{definition}{Definition}[section]

\newtheorem{lemma}[definition]{Lemma}
\newtheorem{proposition}[definition]{Proposition}

\newtheorem{remark}{Remark}
\newtheorem{example}{Example}

\author{F.~Bassetti \footnote{Universit\`a degli Studi di Pavia, Dip.
    Matematica, Via Ferrata 1, 27100, Pavia, Italy;
    e-mail:federico.bassetti@unipv.it }, {M.~Cosentino~Lagomarsino}
  \footnote{
    Universit\`a degli Studi di Milano, Dip.  Fisica and I.N.F.N., 
    Via Celoria 16, 20133
    Milano, Italy; e-mail:marco.cosentino-lagomarsno@unimi.it} { S. Mandr\`a } \footnote{
    Universit\`a degli Studi di Milano, Dip.  Fisica and I.N.F.N., 
    Via Celoria 16, 20133
    Milano, Italy; e-mail:salvatore.mandra@mi.infn.it} }

\title{Exchangeable Random Networks}

\begin{document}


\def \La{ \textbf{\texttt{L}}}
\def \CZ{\mathcal{Z}}
\def \CK{\mathcal{K}}
\def\CF { \mathcal{F}}
\def\CL { \mathcal{L}}
\def\CM { \mathcal{M}}
\def\CH{ \mathcal{H}}
\def\CB {\mathcal{B}}
\def\CP { \mathcal{P}}
\def\CI { \mathcal{I}}
\def\F { \mathcal{F}}
\def\CL { \mathcal{L}}
\def\C { \mathcal{C}}
\def\CW { \mathcal{W}}
\def\CI {\mathcal{I}}
\def\CD {\mathcal{D}}
\def\CN {\mathcal{N}}
\def\CS {\mathcal{S}}
\def\CG {\mathcal{G}}
\def\CP {\mathcal{P}}
\def\CU {\mathcal{U}}
\def\CR {\mathcal{R}}
\def\b {\beta}
\def\a {\alpha}
\def\t {\tau}
\def\d {\delta}
\def\th {\theta}
\def\s {\sigma}
\def\r {\rho}
\def\g {\gamma}
\def\eps {\epsilon}
\def\n {\nu}
\def\a {\alpha}
\def\g {\gamma}
\def\t {\tau}
\def\gg {\eta}
\def\z {\zeta}
\def\lm {\lambda}
\def\I {{\bf I}}
\def\J {\mathbb{I}}
\def\Q {\mathbb{Q}}
\def\N {\mathbb{N}}
\def\G {\mathbb{G}}
\def\E {\mathbb{E}}
\def\RE {\mathbb{R}}
\def\Q {\mathbb{Q}}
\def\N {\mathbb{N}}
\def\Z {\mathbb{Z}}
\def\FF { \mathbb{F}}
\def\q {\diamondsuit}
\def \fine{\diamondsuit}

\def\dirac#1{\J\{#1\}}
\def\parg#1{\left\lbrace#1\right\rbrace}
\def\part#1{\left(#1\right)}
\def\parq#1{\left[#1\right]}
\def\CS{\mathcal{S}}
\def\CG{\mathcal{G}}


\maketitle

\begin{abstract}
  We introduce and study a class of exchangeable random graph
  ensembles. They can be used as statistical null models for empirical
  networks, and as a tool for theoretical investigations. We provide
  general theorems that characterize the degree distribution of the
  ensemble graphs, together with some features that are important for
  applications, such as subgraph distributions and kernel of the
  adjacency matrix.  
  A particular case of directed networks with power-law out--degree is studied
  in more detail, as an example of the flexibility of the model in
  applications.
\end{abstract}

\small \centerline{\bf Key words and phrases} Complex networks, Exchangeable
graphs, Power law distributions, Random Graphs, Scale free graphs,
Transcription networks.

\section{Introduction}

Random graphs have attracted much interest as null- and positive models
for many real-world systems involving many interacting agents, such as
the internet, epidemics, social and biological interactions (see for
instance~\cite{Strogatz01,Newman2003,NewBarabWatts2006,Montoya06}).
In many of these instances, one is naturally confronted with properties
that differ from the classical Erd\"os-R\'enyi model.  We recall that,
in the Erd\"os-R\'enyi model, edges in the graph exist independently
from each other, with a fixed probability (dependent on the dimension of
the graph).  While for the Erd\"os-R\'enyi model analytical expressions
for many of the relevant observable properties of the graph (such as the
diameter, clustering coefficient, component size distributions, subgraph
distribution, giant component, etc) are available, less is known for
other kinds of models.  In the recent years, in connection with the
availability of large-scale data on real-life networks, many studies
addressing random graph models going beyond the Erd\"os--R\'enyi model
have appeared.  Two studies that are worth mentioning are the so-called
"small--world" model~\cite{Watts98} and the preferential-attachment
model~\cite{Barab99}, addressing the empirically observable phenomena of
short shortest-paths and power-law degree distributions respectively.
This new wave of models has affected also the mathematical literature
(see, for instance refs
\cite{Aiello01,Aiello02,Bollobas03,Bollobas03b,Bollobas03c,
Bollobas04,Bollobas04b,Chung03,Chung04,Chung06b,
Newman2003,Newman2003b}).
Among the many recent mathematical books on the subject, we would like
to mention, for classical random graph theory, \cite{Kol98} and
\cite{Bollobas2001}, and, for more recent models of random graphs,
\cite{Chung06} and \cite{durrett}.
From a statistical point of view, which we adopt here, it is natural to
seek a parameterizable stochastic model of complex graphs, that would be
at the same time flexible for practical use and mathematically tractable
for theoretical exploration. Moreover, it is desirable that the
qualitative properties of the model should emerge from some simple
unifying mathematical structure rather than from \emph{ad-hoc}
considerations,
see~\cite{Aiello02,Bollobas2001,Bollobas03b,Bollobas03c,Bollobas04,Chung04,NewBarabWatts2006}.

The aim of this paper is to present a general class of random graphs
that addresses these needs. It was introduced in~\cite{BCBJ} in a
particular case, connected to the study of null models for
transcriptional regulation networks~\cite{BLA+04}. The defining property
of the graph ensemble is the exchangeable structure of its degree
correlations. This symmetry property makes it particularly apt to be
used as a statistical null model.  The most important advantages of such
an approach are the following: (i) Much as in the Erd\"os-R\'enyi model,
some observables can be easily computed analytically for finite sizes
and asymptotically, rather than estimated numerically. (ii) It is fast
and versatile in computational implementations and statistical
applications. As we will show in the different sections of this paper,
many observables that are commonly useful in the analysis of large-scale
networks are particularly simple to access with our ensemble.  In order
to show the range of applicability, we discuss multiple applications to
observables in the model graphs rather than presenting a very detailed
analysis on a single graph feature. In the use as a null model,
differently from other approaches used in the study of transcriptional
and other networks~\cite{IMK+03,RJB96,CDH+05}, our generating method for
random graphs is not designed to conserve the degree sequence of the
observed real graph, but rather as a method to generate graphs with
degree distributions having certain prescribed properties.

The paper is structured as follows. Section~\ref{Sec:1} introduces a rather
general class of random directed network ensembles that can be produced with
the same defining principle of exchangeability, and discusses some simple
variants.  The following part is intended to show how the structure of the
proposed model is useful in the the study of many relevant topological
features of the ensemble. To this aim, in Section~\ref{Sec:2} we prove some
theorems which characterize the degree distributions and the distribution of
the size of the "hub" (or the maximally connected node). In particular, we
show that the model can generate an ensemble characterized by a Poisson limit
distribution for the in--degree, and a mixture of Poisson limit distributions
for the out-degree.  This important property enables to obtain a limit
out--degree distribution with power--law tails.  In the same section, we show
that the probability that the graph is disconnected goes to 1 as the size
of the graph diverges.  Section~\ref{Sec:3} gives some results concerning the
mean number of subgraphs (a quantity of some importance in many applications),
roots and leaves.  Section~\ref{Sec:4} considers a particular Boolean
optimization problem defined on the graph, which emerges in statistical
physics and theoretical computer science.  More precisely, we will give some
results concerning the non-trivial problem of the dimension of the kernel of
the adjacency matrix.  In Section~\ref{Variants} we briefly comment the two
variants of the main model.  Finally, Section~\ref{Sec:7} contains the detailed
analysis of a simple two-parameter ensemble derived from the general model
presented in Section~\ref{Sec:1}.  Some of the proofs are deferred to the
Appendix.

\section{The Model}\label{Sec:1}
\def\A {\mathbb{X}}

Although the ideas we describe are applicable to both directed and
undirected graphs, we will mainly consider here the case of directed
graphs.  Any directed random graph $G_n$ with $n$ nodes is completely
specified by its adjacency matrix $
\A_n=\A(G_n)=[X^{(n)}_{i,j}]_{i,j=1,\dots, n} $, where $X^{(n)}_{i,j}
= 1$ if there is a directed edge $i \rightarrow j$, $0$ otherwise.  
In many applications, such as transcription networks instead of square
matrices, one may also consider rectangular matrices.  The reason for
this is that in some situations it is reasonable to assume that, while
all nodes can receive edges, only a fraction of nodes can send them out
(see \cite{BCBJ} for an introduction to this problem).  Hence, in what
follows we will deal with rectangular matrices $m_n \times n$. As we
will see in Section~\ref{Sec:7}, this is a necessary choice for
networks with power-law degree distributions having exponent equal or
lower than 2 (thus with diverging average) to obtain non--trivial
asymptotics.

One of the interests of our procedure is the fact that it can produce graphs
with different in- and out- degrees distributions.  Naturally, if the graph is
generated by throwing independently each directed edge with a fixed
probability -- as in the case of (undirected) Erd\"os-R\'enyi graphs-- this is
not possible.
In order to build a random graph with different in- and out--degree
distributions, one must give up total independence and allow some kind
of dependence among edges.  In particular, maintaining the maximal
symmetry leads to the choice of exchangeability.

\subsection{Partially Exchangeable Random Graphs}\label{s.sec1.A}

The first general class we will consider includes directed graphs whose in- or
out--degrees, i.e. the columns or the rows of $\A_n$, are exchangeable, while
the out- or in--degrees are stochastically independent.  Differently put, our
model ensemble can be defined using the following generative algorithm. For
each row of $\A_n$, independently, (i) throw a bias $\theta$ from a prescribed
probability distribution $\pi_n$ on $[0,1]$ (ii) and set the row elements of
$\A_n$ to be $0$ or $1$ according to the toss of a coin with bias $\theta$.
Since each row is thrown independently, the resulting probability law is
\begin{equation}
  \begin{split}
    P\{ X_{i,j}^{(n)} &= e_{i,j}, \, i=1,\dots,m_n, j=1,\dots,n \} \\
    &=\prod_{i=1}^{m_n} \int_{[0,1]} \th_i^{\sum_{j=1}^n
      e_{i,j}}(1-\th_i)^{n-\sum_{j=1}^n e_{i,j}} \pi_n(d\th_i) \\
  \end{split}
  \label{eq:P}
\end{equation}
where $e_{i,j} \in \{0,1\}$, $i,j=1,\dots n$.  In other words, each row of $
\A(G_n)$ is independent from the others with exchangeable law directed by
$\pi_n$.  
One can apply an identical procedure to the transposed matrix of
$\A_n$ and switch the role of in- and out--degrees.

It is worth recalling that a random vector, say $(Y_1,\dots,Y_n)$, is
said to be exchangebale if its law is invariant under any permutation,
that is, if for any permutation $\s$ of $\{1,\dots,n\}$,
$(Y_1,\dots,Y_n)$ and $(Y_{\s(1)},\dots,Y_{\s(n)})$ have the same
law. For an introduction to exchangeable sequences and array see,
e.g., \cite{Aldous}.  This hypothesis is important for the use of the
ensemble to produce statistical null models, as it implies symmetry of
the probability distributions with respect to the permutation of
variables, i.e. all the nodes or the agents they represent (genes,
computer routers, etc .) are given an equivalent status.

To complete the model, one has to specify the choice for $\pi_n$, which
determines the behavior of the graph ensemble.
For example, in~\cite{BCBJ}, we have chosen the two-parameter distribution
\begin{equation}\label{defpi}
\pi_n(d\th)=  Z^{-1}_n \th^{-\beta}\J_{(\frac{\a}{n},1]}(\th) d\th
\end{equation}
where $n > \a >0$ and $\beta >1$ are free parameters, $\J_{(\frac{\a}{n},1]}$ is
the indicator function of the interval $(\frac{\a}{n},1]$, taking the
value one inside the interval and zero everywhere else, and
$
Z_n:=((n/\a)^{\beta-1}-1)/(\beta-1)
$
is the normalization constant. As we will see in Section~\ref{Sec:7}, this
choice produces a graph ensemble with heavy-tailed degree sequences.  As a
second example, taking $\pi_n(d\th)=\delta_{\lm/n}(d\th)$, one obtains a
directed version of the Erd\"os-R\'enyi graph. 

A naturally interesting problem is to characterize the general forms of
the probability measure $\pi_n$ that lead to graph ensembles with
qualitatively different characteristics.  In Section~\ref{Sec:2} we
shall give some results in this direction.
Note that a general way of producing the distribution $\pi_n$ for each $n$,
starting form a given "seed" $F$ ($F$ being a fixed distribution function on
$\RE^+$), is easily described by the following assumption:
\begin{equation}\label{semi_parametricmodel}
 F_n(x):= F(x n)/F(n)=\int_{[0,x]} \pi_n(d\th).
\end{equation}
With the above assumption, $F_n$ is a well-defined distribution function on
$[0,1]$ whenever $F(n)>0$, which certainly holds for large enough values of
$n$.

\subsection{Completely exchangeable graphs}\label{s.sec1.B} 

The above described method of generating exchangeable graphs is quite general,
so that one can imagine many simple variants.  For example, one can consider
the following algorithm: (i) throw a bias $\theta$ from a prescribed
probability distribution $\pi_n$ and (ii) set all the elements of $\A_n$ to be
$0$ or $1$ according to the toss of a coin with bias $\theta$.  The resulting
probability law, say $Q$, is
\[
Q\{ X_{i,j}^{(n)} = e_{i,j}; \, i,j=1,\dots,n \}= \int_{[0,1]}
\th^{\sum_{i,j} e_{i,j}}(1-\th)^{n^2-\sum_{i,j} e_{i,j}} \pi_n(d\th)
\qquad
\]
for any $e_{i,j}$ in $\{0,1\}$ $i,j=1,\dots n$, that is
under $Q$, $\{ X_{i,j}^{(n)}; i,j=1,\dots,n \}$
are exchangeable, with de Finetti measure $\pi_n$.

\subsection{Hierarchical models}\label{s.sec1.C} 
Another possible variant considers a hierarchy of probability distributions to
generate the bias of the coins. In this case one can take
\begin{equation}
  \begin{array}{l}
    Q^*\{ X_{i,j}^{(n)} = e_{i,j}, \, i=1,\dots,m_n, j=1,\dots,n \}= \\
     \, \\
    \int_{\RE^+} \prod_{i=1}^{m_n} \int_{[0,1]} \th_i^{\sum_{j=1}^n
      e_{i,j}}(1-\th_i)^{n-\sum_{j=1}^n e_{i,j}} \pi_n(d\th_i|\a) \lambda_n(d\a)
  \end{array}
  \label{eq:Pbis}
\end{equation}
$ \lambda_n$ being a probability on $\RE^+$ and 
$\pi_n(d\th|\a)$ being a kernel on $[0,1]\times \RE^+$, that is: for every $\a$ in $\RE^+$, $\pi_n(\cdot|\a)$ is a 
 measure on the Borel $\s$--field of $[0,1]$ and, for every measurable subset $B$ of $[0,1]$, 
 $\a \mapsto \pi_n(B|\a)$ is measurable. 

\section{Connectivities}\label{Sec:2}

We will carry the main discussion considering the case of partially
exchangeable graphs of Subsection~\ref{s.sec1.B}.  Some brief comments
on the other variants are reported in Section~\ref{Variants}.  In the
rest of the paper, with the exception of Section~\ref{Variants}, we
suppose that all the random elements are defined on the same probability
space $(\Omega,\CF,P)$ and we denote by $\E(Y)$ the mathematical
expectation of a given random variable $Y$ with respect to $P$. With a
slight abuse of notation we shall use indifferently $G_n$, the random
graph, and its adjacency matrix $\A_n=[X^{(n)}_{i,j}]_{i,j}$.

\subsection{In and out connectivity}
The first quantities that we want to characterize are the graph degree
distributions.  The random variable $Z_{m_n,j}:=\sum_{i=1}^{m_n} X^{(n)}_{i,j}
$ represents the in--degree of the $j$-th node in the random graph, while
$S_{n,i}:=\sum_{j=1}^n X^{(n)}_{i,j} $ can be seen as the out--degree of the
$i$-th node ($1 \leq i \leq m_n$).  Note that $(Z_{m_n,1},\dots, Z_{m_n,n})$
are identically distributed as well as $(S_{n,1}, \dots,S_{n,m_n} )$.
Moreover, $(S_{n,1}, \dots,S_{n,m_n} )$ are independent, and each $S_{n,i}$ is
a sum of exchangeable Boolean random variables, while $(Z_{m_n,1},\dots,
Z_{m_n,n})$ are dependent.  
Clearly, the mean degrees are equal to $m_n \mu_n$ and $n \mu_n $,
respectively, where $\mu_n:=P\{X_{i,j}^{(n)}=1\} = \int_{[0,1]} \th \pi_n(d
\th)$ is the probability of the link $i \rightarrow j$.  Note that, while in
the Erd\"os--R\'enyi model $n \mu_n=\lambda$ for every $n$, in this case $n
\mu_n$ generally depends on $n$. On the other hand, when
(\ref{semi_parametricmodel}) is in force, using the well-known fact
that $\E(Y)=\int_0^{+\infty}(1-G(y))dy$
for any positive random variable $Y$ with distribution function $G$, one gets
\[
n \mu_n=\int_{[0,1]} n\th \pi_n(d\th)
 =\int_0^n (1-F(x)/F(n))dx,
\]  
and hence, if $\mu:=\int_0^{+\infty}x dF(x)<+\infty$, it follows that
$n \mu_n=\mu+o(1)$.
The (marginal) degree distributions are given by
\begin{equation}\label{distS_n}
P\{S_{n,i}=k\}= {{n}\choose{k}}\int_{[0,1]} \th^k(1-\th)^{n-k} \pi_n(d
\th) 
\end{equation}
and
\begin{equation}\label{distZ_n}
P\{Z_{m_n,j}=k\}= {{m_n}\choose{k}} \mu_{n}^k(1-\mu_n)^{m_n-k} .
\end{equation}
With the above expressions, the problem of determining the asymptotic
distribution of $(Z_{m_n,1})_{n \geq 1}$ and $(S_{n,1})_{n \geq 1}$ is simply
cast in a central limit problem for triangular arrays. In fact, while for
$(Z_{m_n,1})_{n \geq 1}$ a classical central limit theorem (CLT) for
triangular arrays of independent random variables works, for $(S_{n,i})_{n
  \geq 1}$ one needs a CLT for exchangeable random variables. General CLTs for
exchangeable random variables are well known (see, for instance,
\cite{FLR,RS97}). Here the situation is particularly simple, since we are
dealing with $0-1$ random variables. Consequently, we need only a simple
\emph{ad-hoc} CLT, for exchangeable Boolean random variables.

Let $\tilde\th_n$ be a random variable taking values in $[0,1]$ with
distribution $\pi_n$ and set $T_n:=n \tilde \th_n$. The next proposition shows
that, under a set of reasonable assumptions on $T_n$, the limit law of
$(S_{n,1})_{n \geq 1}$ is a mixture of Poisson distributions, while the limit
law of $(Z_{n,1})_{n \geq 1}$ is a simple Poisson distribution.

\begin{proposition}[CLT]\label{CLT} If $(T_n)_{n \geq 1}$
converges in distribution to a random variable $T$ with distribution
function $F$, then, for every integer $j \geq 1$,
\begin{equation}\label{eq1}
\lim_{n \to +\infty} P\{S_{n,j}=k \} = \E\left[\frac{1}{k!} T^k e^{-T}\right]=
\int_0^{+\infty} \frac{t^k}{k!} e^{-t} dF(t) \qquad (k=0,1,\dots).
\end{equation}
Moreover, if for some $\lambda >0$ and  for a sequence $(a_n)_{n \geq 1}$
\begin{equation}\label{eq1bis}
\lim_{n \to +\infty} a_n \E(T_n)= \lim_{n \to +\infty} n a_n
\int_{[0,1]} \th \pi_n(d\th)= \lambda
\end{equation}
holds true,  then, for every integers $k \geq 0$ and $j$,
\[
\lim_{n \to +\infty} P\{ Z_{m_n,j}=k \} = \frac{\lm^k e^{-\lm}}{k!}
\]
with $m_n=[n a_n]$ {\rm(}$[x]$ being the integer part of $x${\rm)}.
\end{proposition}

\begin{remark}
  (a) If (\ref{semi_parametricmodel}) holds true, then the distribution
    of $T$ is $F$. Indeed, in this case, $$\lim_{n \to
    +\infty}P\{ T_n \leq x \}=\lim_{n \to +\infty}P\{ \tilde \th_n \leq x/n \}
    = \lim_{n \to +\infty} F(x)/F(n)=F(x) \qquad (x \geq 0).$$ (b) It is
    worth noticing that as a corollary of Theorem 5 in \cite{FLR} one
    has that the convergence of $T_n$ is a necessary and sufficient
    condition in order to obtain a Poisson mixture as a limit law for
    $(S_{n,j})_{n \geq 1}$.  Hence, the first part of the previous
    proposition can be proved invoking such a theorem. Nevertheless, for
    the sake of completeness, we shall give here a simple direct proof.
\end{remark}

{\it Proof of Proposition \ref{CLT}.} Since $T_n:=n \tilde \th_n$,
by (\ref{distS_n})
 one
has
\[
P\{S_{n,j}=k \}=\E\left [ {n \choose k} \frac{1}{n^k} T_n^k
  \left(1-\frac{T_n}{n} \right )^{n(1-k/n)} \right ]= \E[\phi_n(T_n)]
\]
where 
\[
\phi_n(x)={n \choose k} \frac{1}{n^k} x^k(1-\frac{x}{n})^{n(1-k/n)}.
\]
Now, 
\[
\E[\phi_n(T_n)]=\E[\phi(T_n)]+R_n
\]
where $\phi(x)=\frac{1}{k!}x^ke^{-x}$ and $R_n=\E[\phi_n(T_n)]-\E[\phi(T_n)]$.
It is plain to check that 
$\phi_n$ converges uniformly on every compact set to $\phi$. Moreover,
since $(T_n)_{n \geq 1}$ converges in distribution, by Prohorov's theorem (see, e.g., Theorem
16.3 in \cite{kallenberg})  
it should be tight, that is for every $\eps>0$ there
exists $K>0$ such that
 $\sup_{n \geq 1} P\{|T_n| \geq K \} \leq \eps $. 
Hence,  one gets that 
\[
\lim_{n \to +\infty}|R_n| \leq \lim_{n \to +\infty } [
\sup_{|x| \leq K} |\phi_n(x)-\phi(x)| + 2 P\{|T_n| \geq K \}] \leq 2\eps. 
\]
At this stage, the first part of the thesis  follows immediately, indeed 
$(T_n)_{n \geq 1}$ converges in distribution if and only if 
$\E[f(T_n)] \to \E[f(T)]$ for every bounded continuous function $f$,
and $\phi$ is bounded and continuous.

The second part of the thesis follows  by the classical Poisson
approximation to binomial distribution  using (\ref{distZ_n}). Indeed
\[
\mu_n= \frac{\lambda}{n a_n}(1+o(1)) 
\]
and $[n a_n]=m_n$ with $n a_n \to +\infty$. To see this last fact, 
observe that, since $T_n$ converges in distribution to $T$,
$\tilde \th_n$ goes to zero  in probability. Using this last fact it is easy to see
that $\E\tilde \th_n=\int_{[0,1]}\th \pi_n(d\th)$ goes to zero, hence  
$n a_n$ must diverge.
$\hfill \fine$

\vskip 0.3cm

\vskip 0.3cm

Since (\ref{eq1}) is a mixture of Poisson distributions with weight
given by $F$, the above result can be used to ``discharge'' the choice
of $\pi_{n}$ on the perhaps more intuitive choice of the mixing
distribution $F$. Clearly, the emergence of heavy-tailed distributions
is not a simple consequence of (\ref{eq:P}), but depends on the choice
of $\pi_n$.  The following example describes a mixing probability which
gives rise to a compact out--degree distribution.

\begin{example}
Take
\[
\pi_n(d\th) =\frac{n\gamma}{1-e^{-\gamma n}} e^{-\gamma n \th}d\th \qquad
(\gamma>0),
\]
or, in other words, assume (\ref{semi_parametricmodel}) with \( F(x)=\int_0^x
\gamma e^{-\gamma t} dt=1-e^{-\gamma x}.  \) With this choice,
according to Proposition \ref{CLT}, the limit
distribution of $S_{n,1}$ is an exponential mixture of Poisson distribution.
Precisely, we find it to be a geometric distribution, indeed
\[
\begin{split}
\lim_{n \to +\infty} P\{S_{n,j}=k \} &= \gamma \int_0^{+\infty}
\frac{t^k}{k!} e^{-t} e^{-\gamma t} dt \\
& =\frac{\gamma}{k!(1+\gamma)^{k+1}} 
 \int_0^{+\infty} y^{(k+1)-1} e^{-y} dy 
 =\frac{\gamma}{k!(1+\gamma)^{k+1}}\Gamma(k+1) \\
&=\frac{\gamma}{1+\gamma}
(1+\gamma)^{-k} \qquad (k=0,1,\dots).\\
\end{split}
\]
Moreover, $a_n=1$ and $\lambda=1/\gamma$ satisfies the 
conditions of  in Proposition \ref{CLT}, yielding
\[
\lim_{n \to +\infty} P\{ Z_{n,1}=k \} =
\frac{\gamma^{-k}
e^{-1/\gamma}}{k!}
\]
As a generalization of the previous example takes, instead of an exponential
distribution, a gamma distribution, i.e.
\[
F(x)=\int_0^x \frac{\gamma^r t^{r-1}e^{-\gamma t}}{\Gamma(r)}dt \qquad (r>0).
\]
It is easy to check that the limiting distribution is a negative binomial
distribution with parameter $r$. That is,
\[
\lim_{n \to +\infty} P\{S_{n,j}=k \} = {r+k-1 \choose k}\left(\frac{\gamma}{1+\gamma}\right)^r
(1+\gamma)^{-k} \qquad (k=0,1,\dots).
\]
Moreover,
\[
\lim_{n \to +\infty} P\{ Z_{n,1}=k \} = \frac{(\frac{r}{\gamma})^{k}
e^{-r/\gamma}}{k!}.
\]

\end{example}

In the above example, mixturing the Poisson distribution with exponential
weights proves insufficient to produce a power-law distribution.  In other
instances, a suitable choice of $F$ in (\ref{eq1}) can give rise to an
out--degree probability distribution with heavy tails. Consider the following


\begin{example}
Assume a slight generalization
of (\ref{defpi}), i.e.
\begin{equation}\label{defpi2}
\pi_n(d\th)=  Z^{-1}_n \th^{-\beta} g(n \th)\J_{(\frac{\a}{n},1]}(\th)
d\th
\end{equation}
with $ 0 < c_1 \leq g(\t) \leq c_2 < +\infty$ for every $\t$ in
$[0,+\infty)$ and
\[
Z_n:= \int_{a/n}^{1} \th^{-\beta} g(n \th) d\th.
\]
Note that (\ref{defpi2}) satisfies 
 (\ref{semi_parametricmodel}) with
 \[
 F(x)=\frac{\int_\a^{x}
{t^{-\b}} g(t) dt}{\int_\a^{+\infty}
t^{-\b}g(t) dt}.
 \]
 Hence, it is straightforward to verify that Proposition \ref{CLT} yields
\[
\lim_{n \to +\infty} P\{S_{n,j}=k \} =\frac{1}{k!}\frac{\int_\a^{+\infty}
{t^{k-\b}} e^{-t}g(t) dt}{\int_\a^{+\infty}
t^{-\b}g(t) dt}=:q_{\a,\b,g}(k).
\]
We now show that such a distribution is a power-law-tailed distribution.
In order to prove this, let us consider first the special case in which
$g=1$, i.e. the older (\ref{defpi}).  With this choice, we get
\[
\lim_{n \to +\infty} P\{S_{n,j}=k \}= \frac{\a^{\b-1}(\b-1)}{k!} \int_{\a}^{+\infty}
t^{k-\b}e^{-t} dt=q_{\a,\b,1}(k)=:p_{\a,\b}(k) \qquad (k \geq 0).
\]
Hence, if $k > \b$, write
\[
p_{\a,\b}(k)=\a^{\b-1}(\b-1)
\left(\frac{\Gamma(k+1-\b)}{\Gamma(k+1)}-\frac{1}{\Gamma(k+1)}\int_0^\a
t^{k-\b}e^{-t} dt \right)
\]
and note that, by the well known asymptotic expansion for the gamma function,
\[
\frac{\Gamma(k+1-\b)}{\Gamma (k+1)} = \frac{1}{k^\b}(1+o(1)) \qquad
\text{as $k \to +\infty$}.
\]
Moreover,
\[
\frac{k^\b}{\Gamma(k+1)}\int_0^\a
t^{k-\b}e^{-t} dt =o(1) \qquad
\text{as $k \to +\infty$}.
\]
Consequently, we get
\[
p_{\a,\b}(k)=\a^{\b-1}(\b-1)\frac{1}{k^\b}(1+o(1)).
\]
Now note that, since
\[
\frac{c_1}{c_2} p_{\a,\b}(k) \leq q_{\a,\b,g}(k) \leq \frac{c_2}{c_1} p_{\a,\b}(k),
\]
$k \mapsto q_{\a,\b,g}(k)$ has power-law tails also for $g \not =1$.
\end{example}

Finally, the following example shows a more complex, already mixtured
distribution, leading to a heavy tail.

\begin{example} Given $\a>1$ and $s>1$, set, for every positive $x$,
\[
f_{\a,s}(x):= \frac{1}{\Gamma(s) \Phi(1,s,\a)}
\int_0^{+\infty} e^{-x(e^{\t}-1)}\t^{s-1}e^{-\t(\a-1)} d\t
\]
where $\Phi(z,s,\a)$ is the well-known Lerch transcendent, defined as
\(
\Phi(z,s,\a):=\sum_{k \geq 0} {z^k}{(\a+k)^{-s}},
\)
for every complex $z$ with $|z|\leq 1$. See, for instance, 9.550 in
\cite{Gradshteyn}.  Note that $f_{\a,s}(x) \geq 0$. Moreover, by means of the
following integral representation 
\(
\Gamma(s)\Phi(z,s,\a)= \int_0^{+\infty}{\t^{s-1}e^{-\t(\a-1)}}(e^{\t}-z)^{-1}d\t
\)
(see 9.556 in \cite{Gradshteyn}), one can check that
\(
\int_0^{+\infty} f_{\a,s}(x) dx=1 \ .
\)
In other words, $f_{\a,s}$ defines a density distribution.  Note that
$f_{\a,s}$ is itself a mixture of exponential densities. Indeed, it can be
rewritten as
\[
\begin{split}
  f_{\a,s}(x) &= \int_0^{+\infty} (e^{\t}-1)e^{-x(e^{\t}-1)}
  \frac{\t^{s-1}e^{-\t (\a-1)}}{\Gamma(s) \Phi(1,s,\a)(e^{\t}-1)} d\t \\
  &= \int_0^{+\infty} u e^{-x u}
  \frac{ \log^{s-1}(u+1)  }{\Gamma(s) \Phi(1,s,\a)(u+1)^{\a+1}} du \\
\end{split}
\] 
with
\[
\int_0^{+\infty} \frac{\t^{s-1}e^{-\t (\a-1)}}{\Gamma(s) \Phi(1,s,\a)(e^{\t}-1)}
d\t = \int_0^{+\infty} \frac{ \log^{s-1}(u+1) }{\Gamma(s)
  \Phi(1,s,\a)(u+1)^{\a+1}} du =1.
\]
It can be verified, with the help of Fubini theorem and 
the already mentioned integral representation of the Lerch
trascendent, that for every real $q$ with $|q|<1$,
\[
\begin{split}
  \sum_{k \geq 0} (iq)^k\int_0^{+\infty} 
\frac{t^{k}}{k!} e^{-t} f_{\a,s}(t)  dt
 &= \int_0^{+\infty} e^{iqt} e^{-t} f_{\a,s}(t) dt \\
  &= \frac{\Phi(iq,s,\a)}{\Phi(1,s,\a)}  = \sum_{k \geq 0} 
\frac{(iq)^k}{\Phi(1,s,\a)(\a+k)^s}\\
\end{split}
\]
(where $i:=\sqrt{-1}$), from which it follows that
\[
\int_0^{+\infty} \frac{t^{k}}{k!} e^{-t} f_{\a,s}(t) dt = \frac{1}{\Phi(1,s,a)(\a+k)^s}.
\]
Hence, if one takes an exchangeable random graph $G_n$, with mixing
distribution satisfying (\ref{semi_parametricmodel}) with
\[
F(x):= \int_0^x f_{\a,s}(t) dt,
\]
then the limit law of $S_{n,1}$ is given by
\[
\lim_{n\to +\infty}P\{S_{n,1}=k\}=\int_0^{+\infty} \frac{t^{k}}{k!} e^{-t}
f_{\a,s}(t) dt=\Phi(1,s,a)^{-1}(\a+k)^{-s}
\]
for every $k \geq 0$.
\end{example}

As the above examples show, the model can produce graphs with disparate
features, depending on the choice of the probability distribution of the coin
biases.
In particular, it is interesting to investigate under which conditions do
heavy-tailed distributions emerge as limit distributions of the out--degree.
If one supposes that $T_n$ converges in law to a random variable with
probability distribution function $F$, we have shown how the question can be
reduced to the problem of determining under which conditions on $F$ the
probability defined by (\ref{eq1}) has heavy tails.  It is worth noticing that
mixtures of Poisson distributions have been extensively studied (see, e.g.,
\cite{Grandell}).  Let us briefly recall some useful properties of such
distributions.  First of all, if
\[
p_k:=\int_0^{+\infty} \frac{1}{k!} t^k e^{-t}dF_i(t) \qquad (k \geq 0,
i=1,2)
\]
for two distribution functions $F_1$ and $F_2$ with $F_i(x)=0$ for every $x
\leq 0$, then $F_1=F_2$, this simple fact was first noticed in
\cite{Feller48}, see also Theorem 2.1 (i) in \cite{Grandell}.  Hence one hopes
to recover many properties of $ p_k:=\int_0^{+\infty} \frac{1}{k!} t^k
e^{-t}dF(t) $ from the properties of $F$.  In particular, Theorem 2.1 in
\cite{Willmot90} states that if $F$ has a density $f$ with respect to the
Lebesgue measure or to the counting measure, such that
\[
 f(x)=L(x) x^{\a} \exp\{-\beta x\}(1+o(1)) \quad \text{as} \,\, x \to +\infty
\]
where $L$ is locally bounded and varies slowly at infinity, 
$\beta \geq 0$, $-\infty <\a <+\infty$ (with $\a<-1$ if $\beta=0$), then
\[
p_k = L(k)\beta^{-(\alpha +1)} (1/(1+\beta))^{k} k^{\a}(1+o(1))
 \qquad \text{as} \,\, k \to + \infty.
\]
Recall that a slowly varying function $L$ is a measurable function
such that
\[
\lim_{x \to +\infty} L(x t)/L(x)=1
\]
for every positive $t$.
Under no assumptions on $F$ we have the following very simple 
\begin{lemma}\label{lemma3.2}
  Let $F$ be a distribution function with $F(x)=0$ for every $x \leq 0$, and
  set
\(
p_k:=\int_0^{+\infty} \frac{1}{k!} t^k e^{-t}dF(t)
\).
Then, for every positive $\gamma$
\[
\sum_{k \geq 0} k^\gamma p_k < +\infty
\]
if and only if
\[
\int_0^{+\infty} t^\gamma dF(t)< +\infty.
\]
\end{lemma}

The proof is deferred to the Appendix.

It is also worth mentioning that a random variable $T$ is a mixture of
Poisson distribution if and only if its generating function
$G_T(s)=\E(s^T)$ is absolutely monotone in $(-\infty,1)$, that is if
$G_T^{(n)}(s) \geq 0$ for every integer $n$ and $s$ in
$(-\infty,1)$. See \cite{Puri} and Proposition 2.2 in ~\cite{Grandell}.
Finally, we recall that the sequence $(p_k)_{k \geq 1}$ inherits many 
properties from $F$.  For
example, $(p_k)_{k \geq 1}$ has a monotone density if $F$ has a monotone
distribution, $(p_k)_{k \geq 1}$ has log-convex density if $F$ has log-convex
distribution, $(p_k)_{k \geq 1}$ is infinite divisible if $F$ is so. For more
details see, for instance, ~\cite{Grandell,Steutel}.

The next subsections will deal with the computation of interesting
observables that go beyond the degree distributions.

\subsection{The hub size}\label{thehub}

As a first example of observable, we discuss the size of the so--called
hub, i.e.  the node having maximal out--degree among the nodes (thus, in
many concrete networks, being the most important for routing and the
most vulnerable to attack, see, e.g. \cite{Bollobas03}).  The hub size
is defined by the expression
\[
H_n:=\max_{i=1,\dots,m_n}(S_{n,i}).
\]
In particular, the most interesting case for the behavior of the hub is when
the tail of the out--degree is power-law, as this means that there can be no
characteristic size for the hub.  As we will explain, it is interesting to give
an analytical expression of the limit law of this quantity under a suitable
rescaling.  The idea is very simple: by stochastic independence, it is clear
that $P\{ H_n \leq x b_n\}=(1- P\{S_{n,1} > x b_n \})^{m_n}$, where $x>0$ is
any positive number.  Now, after setting \( L:=\sup\{ y \geq 0: \limsup_n [y
b_n]/n < 1\}, \) if we can prove that \( P\{S_{n,1} \leq x b_n \}=1- g(x)/m_n
+ o (1/m_n) \) for any $x \leq L$, then
\[
\lim_n P\{ H_n \leq x b_n\} = e^{-g(x)}\J_{[0,L)}(x) + \J_{[L,+\infty)}(x).
\]
We will show that, in some situations, it is possible to determine explicitly
$g$, $b_n$ and $L$.  The following proposition concerns the hub behavior in
case of heavy tails for the out--degree. 

\begin{proposition}\label{prophub}
  Suppose there exist two positive constants $\gg,c_\gg$, a sequence
  of positive numbers $(c_{\gg,n})_{n \geq 1}$, and a sequence of functions
  $(r_n)_{n \geq 1}$, such that, for every $t$ in $(0,1)$
  \begin{equation*}
    \int_{(t,1]}\pi_n(d\th)=c_{\gg,n} \frac{1}{(nt)^\gg} + r_n(t),
  \end{equation*}
$c_{\gg,n} \to c_\gg $ and
  \begin{equation}\label{h2}
    \frac{\int_{0}^1 r_n(t) t^{[b_nx]}(1-t)^{n-[b_nx]-1}dt}
    {B([b_nx]+1,n-[b_nx])}  = o(\frac{1}{m_n}) 
  \end{equation}
with  $b_n:= m_n^{1/\gg}$ and $B(\a,\beta):=\int_{0}^1 u^{\a-1}(1-u)^{\b-1}du$. Then 
\[
\lim_{n \to +\infty} P\{ H_n \leq [ x b_n] \} = e^{- c_\gg x^{-\gg} }\J_{[0,L)}(x) +
\J_{[L,+\infty)}(x)
\]
where
\[
L:=\sup\{ y \geq 0: \limsup_{n \to  +\infty} [y \,m_n^{1/\gg}]/n < 1\}.
\]
\end{proposition}

{\it Proof.} 
First of all let us start recalling the well--known
relation
\begin{equation}\label{eq111}
\sum_{k=[x b_n]+1}^n {n \choose k}  \th^{k}(1-\th)^{n-k}=\frac{\int_0^\th
  t^{[b_nx]}(1-t)^{n-[b_n x]-1}dt}{B([b_nx]+1,n-[b_nx])}
\end{equation}
where $B(\a,\b)=\int_0^1 t^{\a-1}(1-t)^{\b-1}dt
=\Gamma(\a)\Gamma(b)/\Gamma(a+b)$.
{See, e.g., 9.2.5 in \cite{Magnus}.}
Hence, by (\ref{distS_n}), (\ref{eq111}) and Fubini theorem  one gets
\[
\begin{split}
P\{ S_{n,1} > [x b_n] \} & = \int_{[0,1]} \sum_{k=[x b_n]+1}^{n}{n \choose k}
\th^{k}(1-\th)^{n-k}\pi_n(d\th) \\
&=\int_{[0,1]}\int_0^\th  \frac{
  t^{[b_nx]}(1-t)^{n-[b_n x]-1}dt}{B([b_nx]+1,n-[b_nx])} \pi_n(d\th) \\
&
= \frac{1}{B([b_nx]+1,n-[b_nx])}
\int_{0}^1 t^{[b_nx]}(1-t)^{n-[b_n x]-1}F^*_n(t)dt\\
\end{split}
\]
with
\[
F^*_n(t):=\int_{(t,1]}\pi_n(d\th)=P\{\tilde \th_n > t \}.
\]
Now, by hypothesis
\begin{equation*}
F^*_n(t)=c_{\gg,n} \frac{1}{(nt)^\gg} + r_n(t) . 
\end{equation*}
Then
\[
P\{S_{n,1} > [x b_n] \} = \frac{1}{B([b_nx]+1,n-[b_nx])}\frac{c_{n,\gg}}{n^\gg}
\int_0^1 t^{[b_nx]-\gg}(1-t)^{n-[b_n x]-1}dt +R_n(x)
\]
with 
\[
R_n(x):= \frac{1}{B([b_nx]+1,n-[b_nx])} 
\int_{[0,1]} r_n(t) t^{[b_nx]}(1-t)^{n-[b_n x]-1} = o(\frac{1}{m_n}).
\]
Finally, using once more the asymptotic expression $\Gamma(n+a)/\Gamma(n+b)=
n^{a-b}(1+o(1))$ as $n \to +\infty$, one obtains
\[
\begin{split}
P\{S_{n,1} > [x b_n] \} &=\frac{1}{m_n}\left (\frac{c_{n,\gg} m_n}{n^\gg} 
\frac{\Gamma(n+1)\Gamma([b_nx]+1-\gg)}{\Gamma(n+1-\gg)\Gamma([b_nx]+1)}  +o(1)  \right) \\
&=\frac{1}{m_n} \left ( \left(\frac{b_n}{[x b_n]}  \right )^\gg c_\gg  +o(1)
\right) , \\
\end{split}
\]
which is 
\[
P\{S_{n,1} > [x b_n] \}=\frac{1}{m_n} \left[\frac{c_\gg}{x^\gg} +
  o(1)\right].  \qquad \fine
\]

\vskip 0.3cm

We give now two simple conditions that imply the validity of (\ref{h2}), and
can be useful in concrete applications.  The first conditions will be used in
the example that we spell out in detail in the second part of this paper
(Proposition \ref{propHub}).

\begin{lemma}\label{lemmahub1}
  If for some $\a>0$, $C<+\infty$ and $\gg>0$
\begin{equation}\label{h1}
|r_n(t)| \leq C\left (\J\{n t < \a \}(1+(nt)^{-\gg}) + n^{-\gg} \right )
\end{equation}
 then (\ref{h2}) holds true provided that $m_n$ is such that
 $m_n/n^\gg=o(1)$.
\end{lemma}

{\it Proof.} 
Set $\beta_n:=[x b_n]$ and 
\[
I_n(d)= m_n  \frac{1}{B(\b_n+1+d,n-\b_n) } \int_{0}^{\a/n} t^{\b_n+d}(1-t)^{n-\b_n-1} dt.
\]
Hence,
\[
\begin{split}
m_n |R_n(x)| & \leq C m_n \left \{ \frac{1}{B(\b_n+1,n-\b_n) }
\int_{0}^{\a/n} (1+1/(nt)^\gg)t^{\b_n}(1-t)^{n-\b_n-1} dt   +\frac{1}{n^\gg} \right \} \\
&=  C \left \{ I_n(0)+ I_n(-\gg) \frac{ B(\b_n+1-\gg,n-\b_n)}{n^\gg B(\b_n+1,n-\b_n)} +o(1)  \right \} \\
& \leq C \left \{ I_n(0) +  I_n(-\gg)\b_n^{-\gg}+o(1) \right \} .\\
\end{split}
\]
It remains to show that $I_n(0)+ I_n(-\gg)\b_n^{-\gg}=o(1)$. With the help of the Sirling formula, one has
\[
\begin{split}
I_n(d) &= \frac{m_n \Gamma(n+1+d)}{\Gamma(\b_n+1+d)\Gamma(n-\b_n )} \int_0^{\a/n} t^{\b_n+d}(1-t)^{n-\b_n-1}dt \\
&\leq  m_n \frac{m_n \Gamma(n+1+d)}{\Gamma(\b_n+1+d)\Gamma(n-\b_n )} (\a/n)^{\b_n+d-1} \\
&\leq C_1 m_n  \frac{\exp\{ \log(\a)(1+d+\b_n)\} (n+1+d)^{n+d+1/2}}{n^{\b_n+d+1}(\b_n+1+d)^{\b_n+d+1/2}(n-\b_n)^{n-\b_n-1/2}}
\exp\{-(n+1+d)+\b_n+1+d+n-\b_n \}\\
&=C_1 m_n \frac{\exp\{ \log(\a)(1+d+\b_n)\} (1+\frac{1+d}{n})^{n(1+\frac{d+1/2}{n})}}
{\b_n^{\b_n+d+1/2} (1+\frac{1+d}{\b_n})^{\b_n(1+\frac{d+1/2}{b_n})}(1-\frac{b_n}{n})^{n(1-\frac{\b_n+1/2}{n})}}\\
    &\leq C_2 m_n \frac{\exp\{ \log(\a)(1+d+\b_n)-\log(\b_n)(1/2+d+\b_n)\}}{[(1-\frac{\b_n}{n})^{n/\b_n}]^{\b_n-\b_n^2/n} }\\
 &\leq C_3 m_n \exp\{ \log(\a)(1+d+\b_n)-\log(\b_n)(1/2+d+\b_n) +(\b_n -\b_n^2/n) \} \\
 &\leq C_4 \exp\{ \log(m_n)-C_5 \b_n\log(\b_n)  \}. \\
    \end{split}
\]
Since $\b_n= x^{1/\gg} m_n^{1/\gg}(1+o(1))$ and $m_n^{1/\gg}/n=o(1)$ the thesis follows easily. 
$\fine$

\vskip 0.3cm

We conclude this subsection observing that when (\ref{semi_parametricmodel})
is in force, then
\[
\int_{(t,1]} \pi_n(d\th)=1-F_n(t)=\frac{F(n)-1+1-F(nt)}{F(n)},
\]
hence it is natural to assume some hypotheses on $1-F(x)$. In particular,
recall that a distribution function $F$ is in the domain of attraction of the
extreme value Fr\'echet distribution if and only if
\begin{equation}\label{DAfrechet}
\sup \{x: F(x)<1\}=+\infty \qquad  \text{and} 
\,\, 1-F(x)=\frac{1}{x^\gg} L(x) 
\end{equation}
where $L$ is a slowly varying function, see  
 \cite{Galambos}.  This means that
(\ref{DAfrechet}) holds if and only if given a sequence of independent and
identically distributed random variables $(\xi)_{n \geq 1}$ with common law
$F$
\[
\lim_{n \to +\infty }P\{ a_n ^{-1} \max \{\xi_1,\dots, \xi_n\}  \leq x
\}
=e^{- c_\gg x^{-\gg} } \qquad (x >0)
\]
for a suitable normalizing sequence $(a_n)_n$.  In point of fact
(\ref{DAfrechet}) is not sufficient, in our case, to ensure that $r_n$
is a reminder of the right order. Hence, we need a heavier requirement. 

\begin{lemma}\label{lemmahub2} Assume that (\ref{semi_parametricmodel}) is in force
  with
\begin{equation}
1-F(x) =  \frac{c_\gg}{x^\gg} [1+ h(x)]
\end{equation}
for some $\gg>0$ and $0 < c_\gg<+\infty$, 
with 
\[
|h(x)| \leq A \left( \frac{1}{x^{\delta_1}}+ \frac{1}{x^{\delta_2}} \right)
\qquad (x >0, A<+\infty, \delta_1,\delta_2 >0),
\]
then  (\ref{h2}) holds true
with $m_n/n^\gg=o(1)$.
\end{lemma}

{\it Proof.}
Assume that $\delta_1=\delta$ and $\delta_2=\delta$.  In the same
notation of the proof of Proposition \ref{prophub}
\[
r_n(t)=\frac{c_\gg}{F(n)n^\gg}\left( 1+h(n) \right)  
+\frac{c_\gg h(nt)}{F(n)(nt)^\gg}
\]
Now
\[
R_n:=R_n(x)= \frac{1}{B([b_nx]+1,n-[b_nx])} \int_{[0,1]} r_n(t)
t^{[b_nx]}(1-t)^{n-[b_n x]-1}dt = R^{(1)}_n+R^{(2)}_n
\]
with
\[
R^{(1)}_n=\frac{c_\gg (1+h(n))
}{F(n)n^\gg}\frac{1}{B([b_nx]+1,n-[b_nx])}
\int_{[0,1]}t^{[b_nx]}(1-t)^{n-[b_n x]-1}dt= \frac{c_\gg
  (1+h(n))}{F(n)n^\gg}=o(m_n^{-1})
\]
and
\[
R^{(2)}_n=\frac{1}{B([b_nx]+1,n-[b_nx])} \int_{[0,1]}t^{[b_nx]}(1-t)^{n-[b_n
  x]-1} \frac{c_\gg h(nt)}{F(n)(nt)^\gg}dt.
\]
Finally
\[
\begin{split}
  |R^{(2)}_n| & \leq \frac{c_\gg 2A}{n^{\delta+\gg}
    B([b_nx]+1,n-[b_nx])}
  \int_{[0,1]}t^{[b_nx]-(\gg +\delta)}(1-t)^{n-[b_n x]-1} dt \\
  & = c_\gg 2A \frac{ B([b_nx]+1-\gg
    -\delta,n-[b_nx])}{n^{\delta+\gg}  
B([b_nx]+1,n-[b_nx])} \\
  & = c_\gg 2A \frac{ \Gamma([b_nx]+1-\gg -\delta)
    \Gamma(n+1)}{n^{\delta+\gg} \Gamma(n+1-\gg -\delta ) \Gamma(
    [b_nx]+1 )}
  \\
  & \leq c_\gg 2A' \frac{(n+1)^{\delta+\gg} }{n^{\delta+\gg}
    ([b_nx]+1)^{\delta+\gg} } \leq c_\gg 2A' 2^{\delta+\gg} \left
    (\frac{b_n}{[b_nx]}\right)^{\delta+\gg} \frac{1}{b_n^{\delta+\gg}}.
  \\
 \end{split}
\]
for a suitable constant $A'$, hence, since $b_n^{\gg}=m_n$,
$|R^{(2)}_n|=o(m_n^{-1})$. The general case follows in the same way.
$\quad\fine$

\vskip 0.3cm

\begin{example}
As an example it is easy to see that 
\[
1-F(x)=\frac{\alpha^\gg}{(\a+x)^\gg} \qquad (x \geq 0, \a >0, \eta >0)
\]
satisfies the assumption of the previous lemma. 
In point of fact
\[
1-F(x)=\frac{\a^\gg}{x^{\gg}}+ \frac{\a^\gg}{x^{\gg}} \left[ \frac{x^\gg-(\a+x)^\gg}{(x+\a)^\gg}   \right ]
\]
hence
\[
|h(x)| \leq \frac{(\a+x)^\gg-x^\gg}{x^\gg}.
\]
If $x\leq 1$, then  $|h(x)|\leq (1+\a)^\gg/x^\gg$, while if $x>1$
\[
|h(x)| \leq (1+\frac{\a}{x})^\gg-1.
\]
Since $t \to t^\gg$ is a Lipschitz function of constant $\gg(1+\a)^{\gg-1}$ on $[1,1+\a]$,
if $x>1$, it follows that $ (1+\a/x)^\gg-1 \leq \gg(1+\a)^{\gg-1}|1+\a/x-1|=\gg(1+\a)^{\gg-1}a/x$.
Summarizing 
\[
|h(x)|\leq A \left[\frac{1}{x} +\frac{1}{x^\gg} \right ].
\]

\end{example}

\section{{Some non-local features of the graphs}}\label{Sec:3}
In this section, we deal with the subgraphs content and the mean number
of roots and leaves of the model of Subsection \ref{s.sec1.A},

\subsection{Subgraphs}
The simple exchangeable structure of the generated random graphs makes
it possible to compute easily the mean value of the number of subgraphs
"of a given shape" contained in the graph, that can be used for the
discovery of ``network motifs''~\cite{SMM+02,MSI+02,MIK+04,Matiaw06}.

Consider a subgraph, with $k$ nodes and $m$ edge, given by
\[
H=\{ i_1 \to i_{(1,1)}, \dots, i_1 \to i_{(1,m_1)},  i_2 \to i_{(2,1)},
\dots, i_k \to i_{(k,1)}, \dots, i_k \to i_{(k,m_k)} \}
\]
with $\sum_{i=1}^k m_i=m$.
 Of course
\[
P\{H \in G_n\}= \int_{[0,1]} \th_1^{m_1} \pi_n(d\th_1)
\int_{[0,1]} \th_2^{m_2} \pi_n(d \th_2) \dots  \int_{[0,1]}
\th_k^{m_k} \pi_n(d \th_k).
\]
Denote by $\mathfrak{T}$ the set of all  subgraphs isomorphic to $H$
contained in the complete $n$ graph and by  $N(H)$ the cardinality
of such set. Since  the number $\CN_H(G_n)$ of graph isomorphic to $H$
contained in $G_n$, can be clearly written as
\[
\CN_H(G_n)= \sum_{ g \in \mathfrak{T}} \J\{ g \in G_n \},
\]
it follows that
\[
\E[\CN_H(G_n)]= N(H)  P\{H \in G_n\},
\]
indeed by exchangeability $P\{ g \in G_n  \}=P\{H \in G_n\}$ for
every $g$ in  $\mathfrak{T}$.

For example, let us consider the $k$--cycles. A subgraph $H$ is
called $k$--cycle if it has the form
\(
i_1 \to i_2 \to \dots \to i_k \to i_1.
\)
If  $\CN_{C_k}(G_n)$ denote  the number of $k$--cycles contained in
$G_n$, then
\begin{equation}\label{kcili}
\E[\CN_{C_k}(G_n)]=2 {n \choose k}\mu_n^k.
\end{equation}

Things are slightly more complicated for rectangular matrices because in
the evaluation of $N(H)$ one needs to take into consideration also the
constrains given by the fact that only $m_n$ nodes can send outgoing
edges.  In what follows we will discuss mainly the case of square
matrices.

As we shall see, in the study of transcriptional networks, the 3--cycle,
 $ i_1 \to i_2, i_2 \to i_3, i_3 \to i_1 $ is called ``feedback loop''
 (\texttt{fbl}), while, with ``feedforward loop'' (\texttt{ffl}) one
 means a triangle of the form $i_1 \to i_2 \to i_3, i_1 \to i_3$.
 Following the procedure described above, one gets
\begin{equation}\label{fblmean}
\E[\CN_{\mathtt{fbl}}(G_n)]=2 {n
  \choose 3}\mu_n^3.
  \end{equation}
As for the evaluation of feedforward loops, we have
\begin{equation}\label{fflmean}
\E[\CN_\mathtt{ffl}(G_n)] =6 {n \choose 3} \int_{[0,1]} \th^2
\pi_n(d\th)
  \int_{[0,1]} \th \pi_n(d\th).
\end{equation}

It is worth mentioning that in principle it is possible to compute
analytically the variance, as well any other moment of the number of
subgraphs isomorphic to a given subgraph. However, computations become
lengthy and cumbersome rather soon.  As an example, we considered the
variance of the number of feedback loops and feedforward loops.

The key point is evaluating $\E\CN_{\mathtt{ffl}}(G_n)^2$ and $\E
\CN_{\mathtt{fbl}}(G_n)^2$. Again, for the sake of symplicity, we will
deal only with square matrices. It is clear that $\E
\CN_{\mathtt{fbl}}(G_n)^2 =\sum_{t \in \mathfrak{T}} \sum_{s \in
\mathfrak{T}} P\{ s,t \in G_n \} $, $\mathfrak{T}$ being the set of all
feedback loops contained in the complete $n$ graph.  Analogously one
obtains $\E \CN_{\mathtt{fll}}(G_n)^2$ taking as $\mathfrak{T}$ the set
of all feedforward loops. Simple calculations give
\[
\begin{split}
\E[ \CN_{\mathtt{fbl}}(G_n)^2]&=4{n \choose 3}{n-3 \choose 3}\mu_n^6 +
12  {n \choose 3}{n-3 \choose
  2}\mu_n^4
\delta_{2,n} \\
&+6(n-3) {n \choose 3}(\mu_n^3\delta_{2,n}+\mu_n^2\delta_{2,n}^2) +2{n \choose
  3}(\mu_n^3+\delta_2^3) \\
\end{split}
\]
where $\delta_{i,n}:= \int_0^1 \th^i\pi_n(d\th)$. As for $\CN_{\mathtt{ffl}}$,
the computations are longer, but essentially the same. The problem is that
$P\{ s,t \in G_n \}$ can take many different expressions depending on $s$ and
$t$. With straightforward but tedious calculations one gets
\[
\begin{split}
\E[ \CN_{\mathtt{ffl}}(G_n)^2]&=  {n \choose 3} A_n
+ (n-3) {n \choose 3} B_n \\
&+ {n \choose 3}{n-3 \choose 2}C_n +{n \choose 3}{n-3 \choose 3}D_n\\
\end{split}
\]
with
\[
\begin{split}
A_n&=6 \delta_{1,n} \delta_{2,n} + 3\delta_{1,n}^2 \delta_{2,n}
+6 \delta_{1,n}  \delta_{2,n}^2 +3\delta_{2,n}^2 +\delta_{2,n}^3,
 \\
B_n&=  
 30 \delta_{1,n}\delta_{2,n}^2 + 18 \delta_{1,n}^2\delta_{2,n}^2
 +6\delta_{2,n}^3
 +18 \delta_{1,n}^2 \delta_{3,n}
+ 12 \delta_{1,n}\delta_{2,n}\delta_{3,n} + 
 6 \delta_{2,n}\delta_{3,n} + 3 \delta_{3,n}^2, \\
C_n &= 60   \delta_{1,n}^2 \delta_{2,n}^2 
+12      \delta_{2,n}^3 
+24 \delta_{1,n} \delta_{2,n}\delta_{3,n}
+ 12 \delta_{1,n}^2 \delta_{4,n}, \\ 
D_n&=36 \delta_{1,n}^2 \delta_{2,n}^2.\\
\end{split}
\]
Hence,
\begin{equation}\label{varfbl}
\begin{split}
Var(\CN_{\mathtt{fbl}}(G_n))&=12  {n \choose 3}{n-3 \choose
  2}\mu_n^4
\delta_{2,n} \\
&+6(n-3) {n \choose 3}(\mu_n^3\delta_{2,n} +\mu_n^2\delta_{2,n}^2) +2{n \choose
  3}(\mu_n^3+\delta_{2,n}^3)   -4 { n \choose 3}R_n \mu_n^6 \\
\end{split}
\end{equation}
and
\begin{equation}\label{varffl}
\begin{split}
Var(\CN_{\mathtt{ffl}}(G_n))&= {n \choose 3} A_n
+ (n-3) {n \choose 3} B_n \\
&+ {n \choose 3}{n-3 \choose 2}C_n - R_n D_n\\
\end{split}
\end{equation}
with $R_n=[{n \choose 3}-{n -3 \choose 3}]$.

\subsection{Roots and leaves}

We say that $i$ is a root if there is no edge of the kind $j \to i$ but
there is at least one edge of the kind $i \to j$ with $j \not= i$. Loops
do not count. Conversely, we say that $i$ is a leaf if there is no edge
of the kind $i \to j$ but there is at least one edge of the kind $j \to
i$ with $j \not = i$, again we exclude loops and isolated points.  Let
$\mathfrak{L}(G_n)$ be the number of leaves in $G_n$ and
$\mathfrak{R}(G)$ the number of roots in $G_n$. Of course,
$\mathfrak{L}(G_n)=\sum_{i=1}^{n} \mathfrak{L}_{i}(G_n)$ and
$\mathfrak{R}(G_n)=\sum_{i=1}^{m_n} \mathfrak{R}_{i}(G_n)$ where
$\mathfrak{L}_{i}(G_n)$ is equal to $1$ if $i$ is a leaf of $G_n$ and
$0$ otherwise and, similarly, $\mathfrak{R}_{i}(G)=1$ if $i$ is a root
of $G_n$ and $0$ otherwise. It follows that
\[
\mathfrak{R}_{i} (G_n) = \J \left \{ \sum_{j=1}^{m_n} X_{j,i}=0 \right \} \left (1- \J
\left \{ \sum_{j=1, \, i \not = j}^n X_{i,j}=0 \right\} \right),
\]
analogously,
\[
\mathfrak{L}_{i}(G_n)   =  \J \left \{  \sum_{j=1 }^{n} X_{i,j}=0 \right
\}  \left(1-
\J \left \{ \sum_{j=1, j\not=i}^{m_n}  X_{j,i}=0 \right \} \right).
\]
Hence,
\begin{equation}\label{leaves}
\E [\mathfrak{L}_{i}(G_n)  ]= (1-\mu_n)^{m_n} (1-P\{S_{n-1,i}=0 \})
\end{equation}
and
\begin{equation}\label{root}
 \E [\mathfrak{R}_{i}(G_n)  ]= (1-(1-\mu_n)^{m_n-1}) P\{S_{n,i}=0
\}
\end{equation}
and then
\[
\E [\mathfrak{L}(G_n)  ]= n(1-\mu_n)^{m_n} (1-P\{S_{n-1,i}=0 \})
\qquad 
\E [\mathfrak{R}_{i}(G_n)  ]= m_n (1-(1-\mu_n)^{m_n-1}) P\{S_{n,i}=0\}.
\]

\subsection{Connected components}

One of the classic and most studied problems in the mathematics of
random graphs is the existence and the size of the so-called giant
component (see, for instance, refs.~\cite{Bollobas04b,Chung06b,Chung03},
the books \cite{Chung06,durrett} and references therein). 
This is in principle an important property if one wants to use the
ensemble as a null or positive model for a real-world system. In many
empirical instances, such as the Internet, World-Wide Web, and many
biological networks, the existence of a very large component can be
observed directly. For this reason, if this property is absent a model
could have limited applications. 
Of course in our model the existence of a giant
component depends on the choice of the measure $\pi_n$.
A detailed study of this problem is beyond the scope of this work, 
and it  will be dealt with in future papers. 
While, for the moment, we did not prove any general
theorem, in some interesting case, such as the power-law model defined
by (\ref{defpi}), one can study the problem numerically. In this case,
our simulations indicate the emergence of a giant component for all
values of the parameters, which makes the model attractive for
applications (see also \cite{BCBJ}).  On general grounds, it is not
hard to see that for this example the probability that $G_n$ has only
one connected component goes to zero as $n$ diverges (at least for $\b
> 2$ and square matrices). This is a consequence of a more general
proposition.

\begin{proposition}\label{connectedcomp}
Let $m_n=n$ and assume that $\lim_{n \to +\infty }
(1-\mu_n)^{n-1} P\{ S_{n,i}=0 \}=a >0 $, then
\[
\lim_{n \to +\infty}P\{ G_n \,\text{is connected}\, \}=0.
\]
\end{proposition}

{\it Proof.}  
 If $Y(n)=\sum_{i=1}^n Y_{i,n}$ with $Y_{i,n}= \J
\{S_{n,i}=0, Z_{n,i}=0 \}$, then
\[
P\{G_n \,\text{is connected}\, \} \leq P \{ Y(n)=0 \} \leq
\frac{Var(Y(n))}{\E(Y^2(n))}=1-\frac{\E(Y(n))^2}{\E(Y^2(n))}.
\]
Since $\E(Y(n)) =n \E(Y_{1,n})=n P\{ S_{n,i}=0, Z_{n,i}=0 \}$ and
\[
\begin{split}
\E(Y(n)^2)&= n \E(Y_{1,n}^2)
+n(n-1) \E(Y_{1,n}Y_{2,n})\\&=n P\{ S_{n,1}=0, Z_{n,1}=0 \}
+n(n-1) P\{ S_{n,1}=0, Z_{n,1}=0,S_{n,2}=0, Z_{n,2}=0 \}\\
&=n (1-\mu_n)^{n-1}
P\{ S_{n,1}=0\}+n(n-1)(1-\mu_n)^{2n-2}
P\{ S_{n,1}=0\}^2,\\
\end{split}
\]
 we get
\[
\begin{split}
P\{G_n \,\text{is connected}\, \} &\leq
1-\frac{(1-\mu_n)^{2n-2}  P\{ S_{n,1}=0 \}^2}{\frac{n-1}{n}(1-\mu_n)^{2n-2} 
P\{ S_{n,1}=0 \}^2+\frac{1}{n}
(1-\mu_n)^{n-1} P\{S_{n,1}=0 \} }.
\end{split}
\]
Taking the limit for $n \to +\infty$ gives the thesis. $\fine$

\vskip 0.3cm

\section{Threshold properties in the kernel of $A_{n}$ }\label{Sec:4}
Another interesting facet of the exchangeable graph
ensemble is its connection with the theory of systems of random equations
over finite algebraic structures.

This problem has fairly important applications in the theory of finite
state automata, the theory of coding, cryptography and combinatorial
optimization problems (satisfiability, colouring).  This kind of
problems arise in many branches of science, ranging from statistical
physics (theory of glasses) to information theory (e.g. low-density
parity- check codes).  See, e.g.,
\cite{Erdos64,Erdos68,Kol98,MacKay,Mezard02,mezard-2003,Murayama03,Levit05}.

One interesting problem in random linear systems over finite algebraic
structures is to prove a threshold property for the random graph $G_n$
with adjacency matrix $\A_n$ of dimension $m_n \times n$. More
precisely, one aims to prove that if $m_n$ and $n$ diverge with $n/m_n
\to \gamma \leq 1$, then an abrupt change in the behavior of the rank of
the matrix $\A_n$ occurs when the parameter $\gamma$ exceeds a
``critical'' value $\gamma_c$.  This property can be expressed in terms
of the total number of hypercycles in $G_n$ defined as
\begin{equation}\label{eq:CS}
        \CS(\A_n) = 2^{Ker(\A_n)}-1 = 2^{n - m_n}\CN(\A_n)-1
\end{equation}
where $\CN(\A_n)$ is the number of nontrivial (i.e. non zero) solutions of
the linear system in $\G\FF_2$ (the field with elements 0 and 1)
\begin{equation}\label{eq:ls}
        \A^T_n x =_{\G\FF_2} 0.
\end{equation}
Problems of this kind have been extensively studied for a few ensembles
of random graphs, see, for instance, Theorem 3.5.1 in \cite{Kol98} and
Theorem 1 in \cite{Levit05}.
  
In the next proposition, we give an exact expression for the mean value of the
number of solutions of the linear system (\ref{eq:ls}). This expression can be
used to prove the existence of a threshold property for $\CS(\A_n)$.
Moreover, the same expression is a first step for a more exhaustive
characterization of solution space, which shall be dealt with in a forthcoming
paper. All the proofs of this Section are deferred to the Appendix. 

In order to state the next proposition introduce the following notations.
Define
\[ 
        \xi_{n}(i)=\int_{[0,1]}(1-2\th)^i \pi_{n}(d\th) 
\]
and
\[
        \CZ_{n} = \parg{j\in\parg{0,\ 1,\ \dots,\ n}\,|\,\xi_{n}(j) = 0}.
\]

\begin{proposition}\label{1th-moment} 
  Assume that $\A_n$ is a random adjacency matrix of dimension $m_n \times n$
  with law (\ref{eq:P}).
  Then
\begin{equation} 
\E\CN(\A_n) =2^{-n} \sum_{i=1}^{n}{n \choose i}   
(1+\xi_{n}(i))^{m_n} 
\end{equation} 
whenever $\CZ_{n}$ is the empty set.
\end{proposition}

Using the previous result one easily obtains the following large
deviation estimate
 
\begin{proposition}\label{largeDev}
  If $m_n=[\frac{n}{\gamma}]$ $(\gamma \leq 1)$ and $(T_n)_{n \geq 1}$
  converges in distribution to a random variable $T$ with distribution
  function $F$, then
\[ 
\lim_{n \to +\infty} \frac{1}{n} \log (\E\CN(\A_n)) =\sup_{x \in 
  [0,1]}\Theta_\gamma(x)=:I_\gamma
\] 
with
\[ 
\Theta_\gamma(x):= \frac{1}{\gamma}\left [\log\left(1+\int_{[0,+\infty)} e^{-2xt}
dF(t) 
\right)-\gamma\left (x\log(x)+(1-x)\log(1-x)+\log(2) \right)  
  \right ] 
\] 
whenever $\CZ_{n}$ is the empty set for $n$ large enough.
\end{proposition}

Combining the previous result with (\ref{eq:CS}) it is clear that, under
the hypotheses of Proposition \ref{largeDev}, the mean number of
hypercycles $\E\CS(\A_n)$ can be written as
\[
\E \CS(\A_n) =(1+o(1))(
2^{n-m_n}e^{n I_\gamma}P_n-1)=(1+o(1))(\exp\{n(I_\gamma-(1/\gamma-1)\log(2)\}P_n-1)
\]
where $P_n$ is a function of $n$ which is at most polynomial,
i.e. $\frac{1}{n} \log(P_n)=o(1)$ (as $n \to +\infty$).  Hence, if
$I_\gamma > (1/\gamma-1) \log(2)$, it follows that $\E \CS(\A_n)$
diverges exponentially in $n$ as $n$ goes to $+\infty$ , while if
$I_\gamma=(1/\gamma-1)\log(2)$ it is sub--exponential, that is for some
$b \geq 0$, $ \E \CS(\A_n)/n^b$ goes to zero as $n$ diverges.

In point of fact we have the following  
\begin{lemma}\label{cor:thr}
If $\int_{[0,+\infty)}t dF(t)<+\infty$ then 
\begin{equation}\label{C1} 
\sup_{x \in [0,1]}\Theta_\gamma(x) > \Theta_\gamma(0)=
\left(\frac{1}{\gamma}-1\right)\log(2).  
\end{equation} 
If 
\begin{equation}\label{C2} 
\log(x) \left (\int_{[0,+\infty)}t e^{-2xt} dF(t)\right)^{-1}=o(1) \qquad (x \to 0), 
\end{equation} 
then there exists a $\gamma_c$ such that for any $\gamma \leq \gamma_c$
\begin{equation} 
\sup_{x \in [0,1]}\Theta_\gamma(x)= \Theta_\gamma(0)=
\left(\frac{1}{\gamma}-1\right)\log(2),  
\end{equation} 
while for $\gamma > \gamma_c$
\begin{equation}\label{C1} 
\sup_{x \in [0,1]}\Theta_\gamma(x) > \Theta_\gamma(0)=
\left(\frac{1}{\gamma}-1\right)\log(2).  
\end{equation} 
In particular, if 
\begin{equation}\label{C3} 
1-F(t)=t^{-\beta}L(t) 
\end{equation} 
with $0 < \beta <1$ and $L$ a slowly varying function then (\ref{C2}) holds
true.
\end{lemma} 
 
In other words, under the hypotheses of Proposition \ref{largeDev}, if
(\ref{C2}) holds true, then there exists a constant $0<\gamma_c<1$ such that
\[
        \lim_{n\to+\infty} \E\CS(\A_n)/n^b = \left\lbrace
        \begin{array}{ll}
                0 & \text{for some $b=b(\gamma)\geq 0$ if }\,\, \gamma \leq \gamma_c\\
                +\infty & \text{for every $b\geq 0$ if } \,\, \gamma > \gamma_c 
        \end{array}
        \right.
\]
That is, the above mentioned threshold property holds.

\section{Other Models}\label{Variants}
In this short section we give some comments about the other two models
presented in Subsections \ref{s.sec1.B}-\ref{s.sec1.C}.

\subsection{Completely exchangeable graphs}\label{ss4.1} 
Most of the properties and quantities discussed above can be easily
established for the totally exchangeable case.  Again,
$\mu_n:=Q\{X_{i,j}^{(n)}=1\} = \int_{[0,1]} \th \pi_n(d \th)$, and, for
instance, the degree distributions (for a square adjacency matrix), are
given by
\[
Q\{S_{n,i}=k\}= Q\{Z_{n,j}=k\} ={{n}\choose{k}}\int_{[0,1]} \th^k(1-\th)^{n-k}
\pi_n(d \th),
\]
Hence, for instance, we have the following
\begin{proposition} 
  If $(T_n)_{n \geq 1}$ converges in distribution to a random variable $T$
  with distribution function $F$, then, for every integer $j \geq 1$,
\begin{equation*}\label{eq1bisbis}
\lim_{n \to +\infty} Q\{S_{n,j}=k \} = \E\left[\frac{1}{k!} T^k e^{-T}\right]
=\int_0^{+\infty} \frac{t^k}{k!} e^{-t} dF(t) \qquad (k=0,1,\dots) 
\end{equation*}
and 
\begin{equation*}\label{eq1tris}
\lim_{n \to +\infty} Q\{Z_{n,j}=k \} = \E\left[\frac{1}{k!} T^k e^{-T}\right]=\int_0^{+\infty} \frac{t^k}{k!} e^{-t} dF(t) \qquad (k=0,1,\dots). 
\end{equation*}
\end{proposition}

For this model, quantities such as the mean number of subgraphs, roots,
leaves, are again easily computed analytically along the same lines described
above. For example, for motifs
\[
Q\{H \in G_n \}= \int_{[0,1]} \th^{m} \pi(d\th).
\]
when 
\[
H=\{ i_1 \to i_{(1,1)}, \dots, i_1 \to i_{(1,m_1)},  i_2 \to i_{(2,1)},
\dots, i_k \to i_{(k,1)}, \dots, i_k \to i_{(k,m_k)} \}
\]
with $\sum_{i=1}^k m_i=m$.
Hence,
\[
\E_Q(\CN_H(G_n))= N(H)  Q\{H \in G_n\}.
\]
Finally, throwing triangular matrices with the same algorithm, one can easily
generate models for undirected graphs.

\subsection{Hierarchical models}
One interesting use of this variant is that it can be exploited to produce
directed graphs having power-law tailed in- and out-degree distributions with
different exponents. To illustrate this point, we will consider the following
example.
 
\begin{example}
  If $\gamma > \beta >2$, $A>0$,
\(
\lambda_n(d\a) \propto \J_{[A,n/2]} \a^{-\gamma}d\a
\)
and
\(
\pi_n(d\th|\a) \propto \J_{(\a /n,1]} \th^{-\beta} d\th,
\)
then
\(
Q^*\{S_{n,1}=k \} = c_1 k^{-\beta}(1+o(1))
\)
and
\(
Q^*\{Z_{n,1}=k \} = c_2 k^{-\gamma}(1+o(1))
\).
Indeed, it is easy to check (by means of a usual dominated convergence
argument) that
\[
\lim_{k \to + \infty} Q^*\{S_{n,1}=k \}=\frac{(\b-1)(\g-1)A^{\gamma-1}}{k!}
\int_A^{+\infty}\int_\a^{+\infty} {\a^{\b-\g-1}t^{k-\b}e^{-t}} dt d\a
\]
and, moreover,
\[
\frac{(\b-1)(\g-1)A^{\gamma-1}}{k!}
\int_A^{+\infty}\int_\a^{+\infty} {\a^{\b-\g-1}t^{k-\b}e^{-t}} dt d\a=
\frac{\g-1}{\g-\b} p_{A,\b}(k)-\frac{\b-1}{\g-\b} p_{A,\g}(k).
\]
In the same way it is easy to check that
\(
\lim_{k \to + \infty}Q^*\{Z_{n,1}=k \}= p_{u,\g}(k)
\)
with
\(
u=A(\b-1)/(\b-2)
\).
\end{example}

\section{A simple  two parameters model}\label{Sec:7} 

In this section, we focus our attention on random graphs generated by assuming
(\ref{defpi}) and we shall specialize the results of previous sections to this
two parameters model.  This model has been suggested by a biological
application. Hence, before presenting the results, we briefly recall the main
features of a transcription network.

Transcription networks are directed graphs that represent regulatory
interactions between genes.  Specifically, the link $a \rightarrow b$ exists
if the protein coded by gene $a$ affects the transcription of gene $b$ in mRNA
form by binding along DNA in a site upstream of its coding
region~\cite{BLA+04}.
For a few organisms, such as E.~coli and S.~cerevisiae, a significant fraction
of the wiring diagram of this network is
known~\cite{LRR+02,GBB+02,SSG+01,HGL+04}.  The topological features of the
graphs can be studied to infer information on the large-scale architecture and
evolution of gene regulation in living systems.  For instance, the
connectivity and the clustering coefficient have been
considered~\cite{GBB+02}. For this kind of analysis one has to consider null
ensembles of random networks with some topological invariant compared to the
empirical case.
The idea behind it is to establish when and to what extent the empirical
topology deviates from the ``typical case'' statistics of the null ensemble.
For example, a topological feature that has lead to relevant biological
findings, in particular for transcription, is the occurrence of small
subgraphs - or ``network motifs''~\cite{MSI+02,MIK+04,WA03,YSK+04}.

As usual in statistical studies, the choice of the invariant properties for
the randomized counterpart is delicate. For instance, the null ensemble used
to for motif discovery usually conserves the degree sequences of the original
network.  The observed degree sequences for the known transcription networks
roughly follow a power-law distribution for the outdegree, with exponent
between one and two, while being Poissonian in the
indegree~\cite{GBB+02,CLB05}.  These features suggest to consider also
alternative null models for directed random graphs with poisson in degree
distribution and (approximately) power-law out-degree distribution, which can
be easily generated with our model under~\eqref{defpi}. In the remainder of
the paper, we will discuss this case in more detail, showing explicit
calculations of the observables discussed in the previous sections.

\subsection{In and out connectivity}

By simple calculations from ~\eqref{defpi},
we get

\noindent If $1 < \beta < 2$, then
\[
\mu_n=\frac{(\beta-1)\a^{\beta-1}}{(2-\beta)n^{\beta-1}}\frac{
1-\left(\frac{\a}{n} \right)^{2-\beta}}{1-\left (\frac{\a}{n}
\right)^{1-\beta}}
= \frac{1}{n^{\beta-1}}\frac{\a^{\beta-1}(\beta-1)}{2-\beta}
[1+o(1) ] \qquad \text{as} \, n \to +\infty.
\]
If $\beta=2$, then
\[
\mu_n=\frac{\a}{n-\a}(\log n-\log \a) 
= \frac{\log n}{n} \a [1+o(1) ] \qquad \text{as} \, n \to
+\infty.
\]
If $\beta >2$, then
\[
\mu_n=\frac{(\beta-1)}{(\beta-2)}\frac{ \left (\frac{n}{\a}
\right)^{\beta-2}-1 }{\left (\frac{n}{\a} \right)^{\beta-1}-1 }
= \frac{1}{n}\frac{\a(\beta-1)}{\beta-2} [1+o(1) ] \qquad
\text{as} \, n \to +\infty.
\]

The next proposition, which is a consequence of Proposition \ref{CLT}, shows
that $S_{n,i}$ is asymptotically power-law distributed.  While $Z_{n,i}$, at
least with a suitable choice of $m_n$, is asymptotically Poisson distributed.
One has to distinguish among the different possible scalings for $\mu_n$.
More precisely, we have the following

\begin{proposition} Assume that (\ref{defpi}) holds true. Then,
for every $\a >0$ and $\b >1$,
\[
\lim_{n\to +\infty } P\{S_{n,j}=k \}=p_{\a,\b}(k) \qquad (j >0, k
\geq 0).
\]
Moreover, if $\b >2$ and $m_n=[\delta n]$ ($\delta \in (0,1]$, $[y]$ being the
integer part of $y$)
\[
\lim_{n\to +\infty } P\{Z_{m_n,j}=k \}=\frac{e^{-\lm} \lm^k}{k!} \qquad
(j >0, k \geq 0),
\]
where
\(
\lm=\frac{\delta \a (\beta-1)}{(\beta-2)}.
\)
If $\b=2$ and $m_n=[\delta n/\log(n)]$
\[
\lim_{n\to +\infty } P\{Z_{m_n,j}=k \}=\frac{e^{-\delta \a}
(\delta \a)^k}{k!} \qquad (j >0, k \geq 0).
\]
If $1<\b<2$ and $m_n=[\delta n^{\beta-1}]$
\[
\lim_{n\to +\infty } P\{Z_{m_n,j}=k \}=\frac{e^{-\lm}
\lm^k}{k!} \qquad (j >0, k \geq 0)
\]
where
\(
\lm=\frac{\delta \a^{\b-1} (\beta-1)}{(2-\beta)}.
\)
\end{proposition}

It is worth noticing that asking for a degree distribution that brings to an
outdegree having a power-law tail with divergent mean ($\beta \le 2$) poses a
heavy constraint on the number of regulator nodes (the rows of the matrix).

\subsection{Subgraphs}

We will discuss mainly the case of square matrices, where calculations are
simpler and conceptually equivalent.

\subsubsection{$k$--cycles}
Under (\ref{defpi}), using (\ref{kcili}), if $\b >2$
\[
\lim_{n\to +\infty } \frac{1}{2}\E(\CN_{C_k}(G_n)) = \frac{1}{k!}
\left[ \frac{\a (\beta-1)}{(\b-2)} \right]^k,
\]
 if $\b= 2$
\[
\lim_{n\to +\infty }  \frac{1}{2 (\log
n)^k}\E(\CN_{C_k}(G_n))=\frac{\a^k}{k!}
\]
and if  $1<\b<2$
\[
\lim_{n\to +\infty }  \frac{1}{2
  n^{k(2-\b)}}\E(\CN_{C_k}(G_n))=\frac{1}{k!}
\left (\frac{(\b-1)\a^{\b-1}}{2-\b} \right)^k.
\]

\subsubsection{Triangles}
The feedforward loop is a classical example of ``network motif'',
i.e. it is overrepresented in known transcription networks. Conversely,
feedback loops (which in principle could form switches and oscillators)
are usually underrepresented (``anti-motifs'') in transcription
networks~\cite{SMM+02,MIK+04}.

Here, we evaluate, for our model, the mean number of feedback loops versus
feedforward loops.  Under (\ref{defpi}), (\ref{fflmean}) yields
\begin{displaymath}
\E\CN_\mathtt{ffl}(G_n) = 6 {n \choose 3}
\frac{(\b-1)^2}{[(n/\a)^{\b-1}-1]^2 }
\int_{\a/n}^1 \th^{2-\b} d\th \int_{\a/n}^1 \th^{1-\b} d\th.
\end{displaymath}
Hence:

\noindent If $\b >3$, then
\[
\lim_{n\to +\infty } \E(\CN_{\texttt{ffl}}(G_n)) = \frac{(\b-1)^2\a^3}{(\b-3)(\b-2)}
> 3\lim_{n\to +\infty } \E(\CN_{\texttt{fbl}}(G_n)) = \frac{(\b-1)^3\a^3}{(\b-2)^3} .
\]
If $\b= 3$
\[
\lim_{n\to +\infty }  \frac{1}{\log n}\E(\CN_{\texttt{ffl}}(G_n))=\a^3.
\]
If  $2<\b<3$
\[
\lim_{n\to +\infty }
\frac{1}{n^{3-\b}}\E(\CN_{\texttt{ffl}}(G_n))=\frac{\a^\b(\b-1)^2}{(\b-1)(3-\b)}.
\]
If $\b= 2$
\[
\lim_{n\to +\infty }  \frac{1}{n \log n} \E(\CN_{\texttt{ffl}}(G_n))=\a^2.
\]
Finally, if $1<\b< 2$
 \[
\lim_{n\to +\infty }  \frac{1}{n^{5-2 \b} }
\E(\CN_{\texttt{ffl}}(G_n))=\frac{\a^{2\b-2}(\b-1)^2}{(3-\b)(2-\b)}.
\]

At this stage one can give the scaling behavior  of 
ratio of  the mean number of feedback and feedforward loops, which is
\[
\frac{\E\CN_{\mathtt{ffl}}(G_n)}{\E\CN_{\mathtt{fbl}}(G_n)}
 \sim
\left\{%
\begin{array}{ll}
n^{\b-1}  &\,  \, \,  \,\mathrm{if} \ \  1 < \b < 2  \\
n/(\log n)^2 &\,  \, \,  \,\mathrm{if} \ \ \b = 2 \\
    n^{3-\b} &  \,  \, \,  \mathrm{if} \ \  2 < \b < 3 \\
 \log n &  \,  \,  \,\mathrm{if} \ \  \b =3 \\
  \lm &  \,  \, \,\mathrm{if} \ \ \b > 3 \\
\end{array}%
\right.
\]
where $\lm =3(\beta-2)^2(\b-3)^{-1}(\b-1)^{-1} > 1$. Here and in what follows we use $
a_n \sim b_n$ to denote $a_n=b_n(1+o(1))$ as $n\to +\infty$. Thus, the \texttt{ffl}
always dominates, although there is a wide range of regimes. Note that the
dominance of feedforward triangles is even stronger if one considers the
rectangular adjacency matrices discussed above. For example, for $1 < \b < 2$,
and rectangular matrices with $m_n=n^{\b-1}$, we calculate
\[
\frac{ \E\CN_{
    \mathtt{ffl}}(G_n)}{\E\CN_{\mathtt{fbl}}(G_n)} \sim
n.
\]
As for the variances, for instance,
one obtains
\[
Var(\CN_{\mathtt{fbl}}(G_n))
 \sim 
\left\{%
\begin{array}{ll}
n^{5(2-\b)}  &\,  \, \,  \,\mathrm{if} \ \  1 < \b < 2  \\
(\log n)^4 &\,  \, \,  \,\mathrm{if} \ \ \b = 2 \\
\frac{1}{3}(\a\frac{\beta-1}{\beta-2})^3 &  \,  \, \,  \mathrm{if} \ \
\b \geq 2. \\
\end{array}%
\right.
\]

\subsection{Roots and leaves}

By simple computations, from (\ref{root}) we obtain:

\noindent If $1 < \beta < 2$ then
\[
\lim_{n} \frac{1}{n^{2-\b}} \log [ \E( \mathfrak{R}_{i}(G_n) )
]=-\frac{\b-1}{2-\b} \a^{\b-1}
\]
and hence
\(
\E( \mathfrak{R}_{i}(G_n)) \sim   e^{ -\frac{\b-1}{2-\b}
\a^{\b-1} n^{2-\b}}.
\)

\noindent If $\beta=2$ then
\[
\lim_{n} \frac{1}{\log n} \log [ \E( \mathfrak{R}_{i}(G_n) )
]=-\a
\]
 and hence
\(
 \E( \mathfrak{R}_{i}(G_n)) \sim
\frac{1}{n^\a}.
\)

\noindent If $\beta > 2$ then
\[
\lim_{n}  \E( \mathfrak{R}_{i}(G_n) )= (1-e^{-\frac{\b-1}{\b-2}
\a}) p_{\a,\b}(0).
\]

Analogously, from (\ref{leaves}), we derive:

\noindent If $1 < \beta < 2$ then
\[
\lim_{n} \frac{1}{n^{2-\b}} \log [1- \E( \mathfrak{L}_{i}(G_n)
) p_{\a,\b}(0)^{-1} ]=-\frac{\b-1}{2-\b} \a^{\b-1}
\]
and hence
\(
\E( \mathfrak{L}_{i}(G_n)) \sim   (1-e^{ -\frac{\b-1}{2-\b}
\a^{\b-1} n^{2-\b}}) p_{\a,\b}(0).
\)

\noindent If $\beta=2$ then
\[
\lim_{n} \frac{1}{\log n} \log [1- \E( \mathfrak{L}_{i}(G_n)
)p_{\a,\b}(0)^{-1} ]=-\a
\]
 and hence
\(
 \E( \mathfrak{L}_{i}(G_n) )\sim
(1-\frac{1}{n^\a})p_{\a,\b}(0).
\)

\noindent If $\beta > 2$ then
\[
\lim_{n}  \E( \mathfrak{L}_{i}(G_n) )=(1-e^{-\frac{\b-1}{\b-2}
\a})(1- p_{\a,\b}(0)).
\]

Combining all the previous statements, we get
\(
\E(  \mathfrak{L}(G_n) ) \sim n
\)
while
\[
\E(  \mathfrak{R} (G_n) ) \sim
\left\{%
\begin{array}{ll}
   e^{-\lm^2 n^{2-\beta}}  &\, \text{if} \, 1 < \b < 2   \\
    n^{1-\a} &  \, \text{if} \, \b = 2  \\
  n &  \, \text{if} \, \b > 2\\
\end{array}%
\right.
\]
where $\lm^2=\frac{\b-1}{2-\b} \a^{\b-1}$.

In concrete applications, these properties can be used for example to impose a
well-defined scaling for the roots-to-leaves ratio of the null network
ensemble.

\subsection{The hub}

In Section~\ref{thehub}, we have already explored the implications on the
limit laws of the maximally connected node of a power-law distributed
out-degree. Using that results under (\ref{defpi}), it is possible to prove an
explicit limit theorem for the size of the hub.

\begin{proposition}\label{propHub}
For $\b > 2$ and for every positive number $x$
\begin{equation}\label{propHubP1}
\lim_{n \to +\infty}P\{ H_n/b_n  \leq x \} = e^{-(\a/x)^{\b-1}}
\end{equation}
with $m_n=n$ and $b_n=n^{1/(\b-1)}$.
For $\b = 2$ and for every
positive number $x$
\[
\lim_{n \to +\infty}P\{ H_n/b_n  \leq x \} = e^{-(\a/x)^{\b-1}}
\]
with $m_n=b_n=n/\log n$.
Finally, for $1 < \b < 2$ and
$m_n=n^{\b-1}$,
\[
\lim_{n \to +\infty}P\{ H_n/n \leq  x \} =e^{-(\a/x)^{\b-1}}\J_{(0,1)}(x)+\J_{[1,\infty)}(x)
\]
for every positive $x$.
\end{proposition}

\begin{remark}
  {\rm (a) Recall that $e^{-(\a/x)^{\b-1}}\J_{[0,+\infty)}(x)$ is the Frechet
    type II extreme value distribution, that is one of the three kind of
    extreme value distributions which can arise from limit law of maximum of
    independent and identically distributed random variables. 

(b) Note that in the last case the limit distribution is not exactly of
extreme value kind and the probability of finding a hub of size $n$ is
asymptotically finite and equal to $1-e^{-(\a)^{\b-1}}$. This concentration
effect was already noted in~\cite{IMK+03} for another kind of random graphs
ensemble.  }
\end{remark}

{\it Proof of Proposition \ref{propHub}.}
Let $\b>2$.  In the same notation of the proof 
of Proposition \ref{prophub}, 
\[
F^*_n(t)=\J\{\a/n \leq t\}
\frac{\a^{\beta-1}}{n^{\beta-1}-\a^{\beta-1}}(t^{1-\beta}-1)
+\J\{\a/n > t\},
\]
hence $\gg=\b-1$, 
\[
c_{n,\g}= \frac{\a^{\beta-1}}{1-n^{1-\beta}} \to c_\g =\a^{\beta-1}
\]
and
\[
r_n(t)=\J\{\a/n > t\}\left (1-
\frac{1}{(nt)^{\beta-1}}
\frac{\a^{\beta-1}}{1- (\a/n)^{\beta-1}}  \right )
- \J\{\a/n \leq  t\}\frac{1}{n^{\beta-1}}
\frac{\a^{\beta-1}}{1- (\a/n)^{\beta-1}}.
\]
The thesis follows by Proposition \ref{prophub} and Lemma
\ref{lemmahub1} noticing that
\[
|r_n(t)| \leq C \left (
(\frac{1}{(nt)^{\beta-1}}+1)\J\{\a/n > t\} + \frac{1}{n^{\beta-1}} \right ).
\]
Arguing essentially in the same way one can prove the statements for
$\b \leq 2$.  $\quad \fine$

\vskip 0.3cm

For $\b >2$ one can guess that $\E[H_n] \sim n^{1/(\b-1)}$, as claimed
in~\cite{IMK+03} in the analyisis of another scale-free random graph ensemble.
In point of fact, we have the following

\begin{proposition}\label{prophub3} If
$\beta >2$ and $d$ is such that
$\beta-d>1$ then
\[
\lim_{n \to +\infty} \E[n^{-d/(\b-1)}H_n^d ] =(\b-1 )^2\a^2 \Gamma \left( \frac{\b-1-d}{\b-1} \right).
\]
\end{proposition}

{\it Proof.} We begin with the case $d=1$. In Proposition 
\ref{propHub} we have just proved that
$(Y_n)_{n\geq 1}:=(H_n/n^\g)_{n\geq 1}$ converges in distribution 
with $\g=1/(\b-1)$. 
So, it is enough to prove that $(Y_n)_{n\geq 1}$ is uniformly integrable, 
i.e.
\[
\lim_{L \to +\infty} \sup_n \E[|Y_n| \J_{|Y_n| \geq L}]=0.
\]
See for instance Lemma 4.11 in  \cite{kallenberg}.
Note, first, that
\[
\E[|Y_n| \J_{|Y_n| \geq L}] \leq LP\{ H_n/n^\g > L   \} 
+ \int_L^{+\infty}(1-P\{ H_n/n^\g \leq x   \})dx.    
\]
Now by (\ref{propHubP1})
\[
LP\{ H_n/n^\g > L   \} \leq C_1 L(1-e^{-\a^{\b-1}L^{1-\b}}) 
\]
for a suitable constant $C_1$. Hence $\lim_{L \to +\infty} \sup_n LP\{ H_n/n^\g > L   \} =0 $
As for  the second term, setting $F_{S_n}(x)=P\{S_n \leq x \}$, one has 
\[
\begin{split}
1-P\{ H_n/n^\g \leq x   \} &= 1-[F_{S_n}(xn^\g ) ]^n \\
&=1-\exp \{ n \log (F_{S_n}(xn^\g)  )   \} \\
& \qquad [ \text{using $1-e^x \leq -x$}]  \\
&\leq  -n \log \left (1 -(1-F_{S_n}(xn^\g)) \right  ) \\
&\leq (1+C_2) n (1-F_{S_n}(xn^\g)).\\
\end{split}
\]
Hence,
\[
\int_L^{+\infty} (1-P\{ H_n/n^\g \leq x   \}) dx \leq (1+C_2) \int_L^{+\infty} n (1-F_{S_n}(xn^\g)) dx=:I_n.
\]
Since $1-F_{S_n}(xn^\g)=0$ if $xn^\g>n$, that is if $x \geq n^{\frac{\b-2}{\b-1}}$, 
\[
I_n = \frac{(\beta-1)(1+C_2)n}{n^{\b-1}(\frac{1}{\a^{\b-1}}-\frac{1}{n^{\b-1}})} 
\int_L^{n^{\b-2}{\b-1}} \int_{\a/n}^1 \sum_{k=[x n^\g]+1}^n {n \choose k} \th^{k-\b}(1-\th)^{n-k} d\th dx
\]
Now, if $L> (\b-1)/n^\g$ then $[xn^\g]+1>0$ for every $x>L$, hence
\[
\begin{split}
I_n &\leq C_3 n^{2-\b} \int_L^{n^{\b-2}{\b-1}}   \sum_{k=[x n^\g]+1}^n {n \choose k} B(n-k+1,k-\b+1) dx\\
&= C_3 n^{2-\b}
\int_L^{n^{\b-2}{\b-1}}  \frac{\Gamma(n+1)}{\Gamma(n-\b+2)} \sum_{k=[x n^\g]+1}^n  
\frac{\Gamma(k-\beta+1)}{\Gamma(k+1) }dx \\
& \leq C_4 n \int_L^{n^{\b-2}{\b-1}} \sum_{k=[x n^\g]+1}^n  \frac{1}{k^\b} dx \\
\end{split}
\]
at least for $L$ large enough. 
Since,
\[
\sum_{k=M}^n \frac{1}{k^\b} \leq \int_{M-1}^n \frac{1}{x^\b} dx 
\]
it follows that
\[
I_n \leq C_4 \int_L^{+\infty} \left ( \frac{1}{[xn^\g]} \right)dx 
\leq C_5 \frac{1}{(L-1)^{\beta-2}}.
\]
The proof of the case with $d>1$ follows an identical procedure, with
$x^{1/d}$ in place of $x$ and $L^{1/d}$ in place of $L$.  $\fine$

\vskip 0.3cm

\subsection{Random linear system in $\G\FF_2$}
Under  (\ref{defpi}) one has 
\[
F(x)= \a^{\b-1 }(\b-1) \int_\a^x \frac{1}{t^{\b}} dt=   
 \left(1-  \frac{\a^{\b-1 }}{x^{\b-1 }}\right) \qquad (x>\a). 
\]
Hence, applying Lemma \ref{cor:thr}, one has that, if $1<\b<2$, then there exists a constant $\gamma_c(\beta)$ 
such that 
\[
        \lim_{n\to+\infty} \E\CS(\A_n)/n^b = \left\lbrace
        \begin{array}{ll}
                0 & \text{for some $b=b(\gamma)$ if }\,\, \gamma \leq \gamma_c(\beta)\\
                +\infty & \text{for every $b\geq 0$ if } \,\, \gamma > \gamma_c(\beta) 
        \end{array}
        \right.
\]
While if $\b> 2$ no threshold property holds since $\int_0^{+\infty}x dF(x)<+\infty$.

\appendix
\section{Appendix}

{\it Proof of Lemma \ref{lemma3.2}.}  Let $k>\gamma$,
$k$ being an integer.  By hypothesis
\[
\frac{\Gamma(k+1)}{\Gamma{(k-\gamma+1)}}p_k=\int_0^{+\infty} 
\frac{t^{k-\gamma}}{\Gamma(k-\gamma+1)} e^{-t}dG(t)
\]
where $G(x)=\int_{(0,x]} t^{\gamma} dF(t)$. Summing both sides on $k$, one can
write
\[
\sum_{k=[\gamma]+1}^{M+[\gamma]+1} \frac{\Gamma(k+1)}{\Gamma(k-\gamma+1)}p_k= 
\int_0^{+\infty} \phi_{\gamma,M}(t) e^{-t} dG(t)
\]
with
\[
\phi_{\gamma,M}(t):=\sum_{k=[\gamma]+1}^{M+[\gamma]+1} 
\frac{t^{k-\gamma+1}}{\Gamma(k-\gamma+1)} =\sum_{m=0}^{M}
\frac{ t^{m+\nu_\gamma}}{\Gamma(m+\nu_\gamma+1)},
\]
and $\nu_\gamma:=[\gamma]+1-\gamma$.
Hence,
\begin{equation}\label{ee1}
\sum_{m=0}^{M}
\frac{\Gamma(m+[\gamma]+2)}{\Gamma(m+\nu+1)}p_{m+1+[\gamma]}= \int_0^{+\infty}
\phi_{\gamma,M}(t) e^{-t} dG(t).
\end{equation}
Now, for every $t>0$ and $\nu_\gamma$ in $(0,1)$, by 5.2.7.20 in \cite{Prudnikov}, one has
\[
\lim_{M \to +\infty}\phi_{\gamma,M}(t)=\sum_{m=0}^{+\infty}
\frac{ t^{m+\nu_\gamma}}{\Gamma(m+\nu_\gamma+1)}=g(\nu_\gamma,t)e^t
\]
where
\[
g(\nu_\gamma,x)=\frac{1}{\Gamma(\nu_\gamma)} \int_0^x \t^{\nu_\gamma-1}e^{-\t} d\t.
\]
Moreover, $\phi_{\gamma,M}(t) \geq 0$ and the convergence is clearly monotone.
Hence, taking the limit as $M$ goes to
$+\infty$ in (\ref{ee1}), by monotone convergence one obtains
\[
\sum_{m=0}^{+\infty}\frac{\Gamma(m+[\gamma]+2)}
{\Gamma(m+\nu_\g+1)}p_{m+1+[\gamma]} =\int_0^{+\infty} g(\nu_\gamma,t) t^\gamma
dF(t),
\]
with $\int_0^{+\infty} g(\nu_\gamma,t) t^\gamma dF(t)<+\infty$ if and only if
$\sum_{m=0}^{+\infty}{\Gamma(m+[\gamma]+2)}{\Gamma(m+\nu_\g+1)^{-1}}
p_{m+1+[\gamma]}<+\infty$.  
Now, since $g(\nu_\gamma,x)$ is a distribution function, one has $\int_0^{+\infty}
g(\nu_\gamma,t) t^\gamma dF(t)<+\infty$ if and only if $\int_0^{+\infty} t^\gamma
dF(t)<+\infty$. Moreover, since 
${\Gamma(m+[\gamma]+2)}{\Gamma(m+\nu_\g+1)^{-1}}=(m+[\gamma]+1)^\gamma p_{m+1+[\gamma]}(1+o(1))$
as $m \to +\infty$, $\sum_{m=0}^{+\infty}
{\Gamma(m+[\gamma]+2)}{\Gamma(m+\nu_\g+1)^{-1}}p_{m+1+[\gamma]} < +\infty$ if and
only if $\sum_{m=0}^{+\infty} (m+[\gamma]+1)^\gamma p_{m+1+[\gamma]}<+\infty$,
which proves the lemma.  
$\quad \fine$

{\it Proof of Proposition \ref{1th-moment}.} Denote by $M(m_n,\,n)$ the set of all $m_n\times n$--adjacency matrices.  
The number of solutions of linear system $\A^T_n x =_{\G\FF_2} 0$ is defined as
\[
       \CN(\A_n) = 
        \sum_{x\in\G\FF_2^{m_n}}\dirac{\A_n^T x=_{\G\FF_2} 0}.
\]
Now note that
\[        
        \dirac{x =_{\G\FF_2} 0} = \frac{1+(-1)^x}{2}
\]
and write 
\[ 
        \E\CN(\A_n) = 
                \sum_{A_n \in M(m_n, n)} P\{\A_n= A_n\}
                \sum_{x\in\G\FF_2^{m_n}} \prod_{j=1}^{n}\frac{1+(-1)^{\sum_{i=1}^{m_n} (A_n)_{ij}x_i}}{2}. 
\] 
Using (\ref{eq:P}) rewrite the last expression as 
\begin{align*}
        \E&\CN(\A_n) \\ 
        &= 2^{-n} \sum_{x\in\G\FF_2^{m_n}}\sum_{A_n \in M(m_n, n)}
                \int_{[0,1]^{m_n}}\parq{\prod_{i=1}^{m_n} \pi_{n}(d\theta_i)} 
                \parq{\prod_{j=1}^{n}\part{1+(-1)^{\sum_{i=1}^{m_n} (A_n)_{ij}x_i}}} 
                 \parq{\prod_{j=1}^{n}\prod_{i=1}^{m_n}\theta_i^{(A_n)_{ij}}\part{1-\theta_i}^{1-(A_n)_{ij}}}\\ 
            &= 2^{-n} \sum_{x\in\G\FF_2^{m_n}}\int_{[0,1]^{m_n}}\prod_{i=1}^{m_n}
\pi_{n}(d\theta_i) 
                \parq{\prod_{j=1}^{n}\sum_{(A_n)_j\in\G\FF_2^{m_n}}
                \part{1+(-1)^{\sum_{i=1}^{m_n} (A_n)_{ij}x_i}} 
                \prod_{i=1}^{m_n}\theta_i^{(A_n)_{ij}}\part{1-\theta_i}^{1-(A_n)_{ij}}} 
\end{align*} 
where $(A_n)_j = \{(A_n)_{1j},\ \dots, (A_n)_{m_n j}\}$
and $(A_n)_{ij}$ is the element in position $(i, j)$ of matrix $A_n$. 
Since the above expression in square brackets  is independent of $j$, $\E\CN(\A_n)$ can be written as
\begin{align*} 
        \E\CN&(\A_n) \\ 
        &= 2^{-n} \sum_{x\in\G\FF_2^{m_n}}\int_{[0,1]^{m_n}}\parq{\prod_{i=1}^{m_n}
\pi_{n}(d\theta_i)} 
                \parq{\sum_{a\in\G\FF_2^{m_n}}\part{1+(-1)^{\sum_{i=1}^{m_n} a_{i}x_i}} 
                \prod_{i=1}^{m_n}\theta_i^{a_{i}}\part{1-\theta_i}^{1-a_{i}}}^{n}.
\end{align*}
At this stage note that 
\[
        \sum_{a\in\G\FF_2^{m_n}}\prod_{i=1}^{m_n}\theta_i^{a_{i}}\part{1-\theta_i}^{1-a_{i}}
= 1
\]
and then
\begin{align*} 
        \E&\CN(\A_n) = 2^{-n} \sum_{x\in\G\FF_2^{m_n}}\int_{[0,1]^{m_n}}\parq{\prod_{i=1}^{m_n}
\pi_{n}(d\theta_i)} 
                \parq{1+\sum_{a\in\G\FF_2^{m_n}}(-1)^{\sum_{i=1}^{m_n} a_{i}x_i} 
                \prod_{i=1}^{m_n}\theta_i^{a_{i}}\part{1-\theta_i}^{1-a_{i}}}^{n}\\ 
        &= 2^{-n} \sum_{x\in\G\FF_2^{m_n}}\int_{[0,1]^{m_n}}\parq{\prod_{i=1}^{m_n}
\pi_{n}(d\theta_i)} 
                \parq{1+\prod_{i=1}^{m_n}\sum_{a_i\in\G\FF_2}\left((-1)^{x_i} 
                \theta_i\right)^{a_{i}}\part{1-\theta_i}^{1-a_{i}}}^{n}
\end{align*}
and after summing over $a_i$ we have
\begin{align}\label{eq:vm1}
        \E\CN(\A_n) 
        = 2^{-n} \sum_{x\in\G\FF_2^{m_n}}\int_{[0,1]^{m_n}}\parq{\prod_{i=1}^{m_n}
\pi_{n}(d\theta_i)} 
                \parq{1+\prod_{i=1}^{m_n} \part{1-\theta_i\part{1-(-1)^{x_i}}}}^{n}.
\end{align}
Now, using
\[
        \dirac{\bar x =_{\G\FF_2} 0} = \frac{1-(-1)^x}{2}
\]
where $\bar x = x + 1$ in $\G\FF_2$, expression (\ref{eq:vm1}) can be written as
\begin{align*}
        \E\CN(\A_n)  
        = 2^{-n} \sum_{x\in\G\FF_2^{m_n}}\int_{[0,1]^{m_n}}\parq{\prod_{i=1}^{m_n}
\pi_{n}(d\theta_i)} 
                \parq{1+\prod_{i=1}^{m_n}\part{1-2\theta_i\,\dirac{\bar x_i =_{\G\FF_2} 0}}}^{n}.
\end{align*} 
Moreover, since
\[
        \part{1-2\theta_i\,\dirac{\bar x_i =_{\G\FF_2} 0}} = 
        \part{1-2\theta_i}^{\dirac{\bar x_i =_{\G\FF_2} = 0}}
\]
we can rewrite the mean number as 
\begin{align*}
        \E\CN(\A_n)  
        = 2^{-n} \sum_{x\in\G\FF_2^{m_n}}\int_{[0,1]^{m_n}}\parq{\prod_{i=1}^{m_n}
\pi_{n}(d\theta_i)} 
                \parq{1+\prod_{i=1}^{m_n}\part{1-2\theta_i}^{\dirac{\bar x_i =_{\G\FF_2} 0}}}^{n}.
\end{align*}
After the expansion of the last square bracket we obtain 
\begin{align*} 
        \E\CN(\A_n) &=  
                2^{-n} \sum_{j=1}^{n}\binom{n}{j}\prod_{i=1}^{m_n} \sum_{x_i\in\G\FF_2}
\int_{[0,1]} \pi_{n}(d\theta_i) 
                \part{1-2\theta_i}^{\dirac{\bar x_i =_{\G\FF_2} 0} j}\\ 
        &= 2^{-n} \sum_{j=1}^{n}\binom{n}{j}\prod_{i=1}^{m_n} \sum_{x_i \in \G\FF_2} \xi_{n} 
                \part{\dirac{\bar x_i =_{\G\FF_2} 0} j}. 
\end{align*} 
Finally, it is easy to see that the last sum is independent of $i$. Then 
\begin{align*} 
        \E\CN(\A_n) = 
                2^{-n} \sum_{j=1}^{n}\binom{n}{j}\parq{\sum_{\sigma\in\G\FF_2} \xi_{n} 
                \part{\dirac{\bar \sigma =_{\G\FF_2} 0} j}}^{m_n} 
                &= 2^{-n} \sum_{j=1}^{n}\binom{n}{j}\part{1 + \xi_{n}(j)}^{m_n}.
\end{align*} 
$\fine$

{\it Proof Proposition \ref{largeDev}.} 
First of all, observe that 
\[ 
\E\CN(\A_n)= \sum_{j=1}^{n}  2^{-n}{n \choose j} \exp\{ n \psi_{n}(j/n)\} 
\] 
where 
\[ 
\psi_{n}(x)=\frac{m_n}{n}\log(1+\xi_{n}(x n))=\frac{m_n}{n} 
\log(1+\E[(1-\frac{2T_{n}}{n})^{x n}]). 
\] 
Now recall that one of the most classical example  
of large deviation estimate is   
\[ 
\lim_{M \to +\infty } 
\frac{1}{M} \log \left (  \sum_{j=0}^M {M \choose j} e^{M f_{M}(j/M) }  \right )  
=\sup_{x \in [0,1]} [f(x) -\{ x \log (x) +(1-x)\log (1-x)+\log(2) \}  ] 
\] 
whenever 
\( 
\lim_{M \to +\infty } \sup_{ x \in [0,1]} | f_{M}(x)-f(x) |=0 
\), $f$ being a continuous function on $[0,1]$. See, e.g., Theorem 7.1 and 10.2 in \cite{ellislecture}.  
Hence, the thesis follows if we prove that  for every $K<+\infty$ 
\begin{equation}\label{eq:M5} 
\lim_{n \to + \infty} \sup_{|x| \leq K } 
|\psi_{n}(x)-\frac{1}{\gamma}\log\left(1+\int_{[0,+\infty)} e^{-2xt} dF(t)\right)|=0. 
\end{equation} 
To prove (\ref{eq:M5}) it is enough to  
 prove that 
for every $K<+\infty$ 
\begin{equation}\label{eq:M1} 
\lim_{M \to + \infty}  
\sup_{|x| \leq K }  |\E[(1-\frac{2T_M}{M})^{Mx}- e^{-2Tx}] |  =0. 
\end{equation} 
Since $T_M$ converges weakly to $T$ and $e^{-t}$ is a bounded and continuous   
function on  $[0,+\infty)$, then 
\begin{equation}\label{eq:M3} 
\lim_{M \to +\infty} \E| e^{-2T_M} -e^{-2T}| = 0. 
\end{equation} 
Moreover we claim that 
\begin{equation}\label{eq:M2} 
\lim_{M \to +\infty}\E|(1-\frac{2T_M}{M})^M- e^{-2T_M} | = 0. 
\end{equation} 
To prove this last claim, set 
\( 
\phi_n(x)=(1-\frac{x}{n})^{n} 
\) 
and note that  $\phi_n$ converges uniformly on every compact set to $e^{-x}$. 
Hence, given $K$,   
\[ 
\lim_{M \to +\infty} \sup_{|x| \leq K} |\phi_M(x)-e^{-x}| =0. 
\] 
Moreover, since $(T_M)_{M \geq 1}$ is tight,  for every $\eps$ there exists  
$K>0$ such that  
\( 
\sup_{M \geq 1} P\{|T_M| \geq K \} \leq  {\eps} 
\). 
Now   
\( 
|(1-\frac{2T_M}{M})^M- e^{-2T_M} | \leq 2 
\) 
and then  
\[ 
\lim_{M \to +\infty}  \E |(1-\frac{2T_M}{M})^M- e^{-2T_M} | 
\leq \lim_{M \to +\infty } [ 
\sup_{|x| \leq K} |\phi_M(x)-\phi(x)| + 2 P\{|T_M| \geq K \}] \leq 2\eps. 
\] 
That is (\ref{eq:M2}). 
Finally, given $a,b$ in $[-1,1]$ and $x>0$  
\[ 
|a^x-b^x| \leq \sup_{y \in [-1,1]} |\frac{d}{dy} y^x |  |a-b| = x |a-b|, 
\] 
hence, since $0< T_M \leq M$, one has 
\( 
-1\leq 1-\frac{2T_M}{M} \leq 1 
\) 
and then 
\begin{equation}\label{eq:M4} 
\begin{split} 
|\E[(1-\frac{2T_M}{M})^{Mx}- e^{-2Tx}] |   
&\leq \E  
|( 1-\frac{2T_M}{M})^{Mx}- e^{-2T_M x} | + \E| e^{-2T_Mx} -e^{-2Tx}| \\ 
& \leq |x| \{\E|(1-\frac{2T_M}{M})^M- e^{-2T_M} | + \E| e^{-2T_M} -e^{-2T}| \} .\\ 
\end{split} 
\end{equation} 
Combining (\ref{eq:M3}), (\ref{eq:M2}) and (\ref{eq:M4}) we get 
 (\ref{eq:M1}).  
$\fine$ 
\\

{\it Proof of Lemma \ref{cor:thr}.} 
Note that,
for every $x$ in $(0,1)$,  
\[ 
\gamma\, \frac{d}{dx} \Theta_\gamma(x)=\frac{2\int_{[0,+\infty)}t e^{-2xt}
dF(t)}{1+\int_{[0,+\infty)} e^{-2xt} dF(t)}-\gamma\log x+\gamma \log(1-x). 
\] 
Hence, if $\int_{[0,+\infty)}t dF(t)<+\infty$ then $\lim_{x \to 0^+}  \frac{d}{dx}
\Theta(x)=+\infty$ and then, 
$\Theta$ is strictly increasing in a neighborhood of $0$. This last fact implies that 
\[ 
\sup_{x \in [0,1]}\Theta_\gamma(x) > \Theta_\gamma(0)=
\left(\frac{1}{\gamma}-1\right)\log(2).  
\] 
If (\ref{C2}) holds true then 
 $\lim_{x \to 0^+}  \frac{d}{dx} \Theta_\gamma(x)=-\infty$, and hence, there exists
$\gamma_c$ such that for any $\gamma \leq \gamma_c$  
\[ 
\sup_{x \in [0,1]}\Theta_\gamma(x)= \Theta_\gamma(0)=
\left(\frac{1}{\gamma}-1\right)\log(2).  
\] 
Now set 
\[ A(x)=\int_0^x t dF(t),  \qquad \text{and} \qquad 
H(s):= \int_0^{+\infty} t e^{-ts} dF(t)=\int_0^{+\infty} e^{-ts} dA(t). 
\] 
The well--known Karamata tauberian theorem (see, e.g. \cite{FellerVolII}) 
 yields that, given $\sigma >0$ and 
$L$ slowly varying, 
\( 
H(s) \sim s^{-\sigma}L(1/s)  
\) 
as $s$ goes to $0$ if and only if 
$A(x) \sim x^{\sigma}L(x)/\Gamma(1+\sigma)$ as $x $ goes to $+\infty$.  
Hence, it remains to prove that if (\ref{C3}) holds true then 
$A(x) \sim x^{\sigma}L(x)/\Gamma(1+\sigma)$. 
Observe that 
\[ 
A(x)=-L(x)x^{1-\beta}+ \int_0^x s^{-\beta}L(s)ds. 
\] 
At this stage the claim follows since  
 it is easy to check that 
 $\int_0^x s^{-\beta}L(s)ds=x^{1-\beta}\tilde L(x)$, where  
$\tilde L(x)$ is still slowly varying. $\quad \fine$      
\\

\section*{Acknowledgments}
We would like to thank Bruno Bassetti for useful discussions and for
having encouraged us during this work. We are also grateful to the referees for many helpful comments.


\end{document}